# Snake locomotion learning search

Sheng-Xue He

(Business School, University of Shanghai for Science and Technology, Shanghai 200093, China)

**Abstract:** This research introduces a novel heuristic algorithm known as the Snake Locomotion Learning Search algorithm (SLLS) designed to address optimization problems. The SLLS draws inspiration from the locomotion patterns observed in snakes, particularly serpentine and caterpillar locomotion. We leverage these two modes of snake locomotion to devise two distinct search mechanisms within the SLLS. The first mechanism, inspired by serpentine locomotion, primarily serves the purpose of exploration within the new algorithm. Unlike conventional methods with explicit destinations, this approach adopts serpentine locomotion's agility to explore the search space effectively. In contrast, the second mechanism draws inspiration from caterpillar locomotion and focuses on exploiting the search space while utilizing historical information during the optimization process. To emulate the way that snakes use pheromones for information exchange, we implement an ordered list of visible spots to record historically successful solutions. This list evolves over time, and snakes employ it to guide their caterpillar-like locomotion. In our quest to mimic a snake's natural adaptation to its surroundings, we incorporate a learning efficiency component generated from the Sigmoid function. This helps strike a balance between exploration and exploitation capabilities throughout the SLLS computation process. The efficacy and effectiveness of this innovative algorithm are demonstrated through its application to 60 standard benchmark optimization problems and seven well-known engineering optimization problems. The performance analysis reveals that in most cases, the SLLS outperforms other algorithms, and even in the remaining scenarios, it exhibits robust performance. This conforms to the No Free Lunch Theorem, affirming that the SLLS stands as a valuable heuristic algorithm with significant potential for effectively addressing specific optimization challenges.

**Key words:** heuristic algorithm; snake locomotion; learning algorithm; engineering optimization; swarm intelligence

## 1. Introduction

Complex optimization problems frequently arise in real-world scenarios, often assuming intricate forms that pose formidable challenges for classical mathematical optimization techniques[1]. These traditional methods typically rely on derivative information or specific features of functions. However, when such information is unavailable or insufficient, heuristic algorithms emerge as a valuable solution. Heuristic algorithms operate as versatile "black-box" tools, requiring minimal or no prior knowledge about the optimization problem. Their successful application across diverse fields over the past few decades underscores their effectiveness in addressing these intricate challenges.

Nevertheless, as elucidated by the "No Free Lunch Theorems for optimization" [2], no single general-purpose algorithm can excel in all cases. This theorem explicitly demonstrates that an algorithm's performance gain in one problem class inevitably entails a trade-off in performance on others. This theoretical insight elucidates the persistent demand for novel heuristic algorithms in practice.

In this paper, we introduce a novel heuristic algorithm known as the Snake Locomotion Learning Search algorithm (SLLS). Inspired by the search behaviors of snakes within their habitats, the SLLS employs two key modes of snake locomotion to enhance its exploration and exploitation capacities. Specifically, it draws inspiration from serpentine locomotion for agile exploration without requiring explicit destinations. Additionally, caterpillar locomotion is harnessed to exploit the search space while leveraging historical information during the optimization process.

The SLLS also incorporates a communication mechanism reminiscent of pheromone exchange among snakes. This mechanism serves to store historical global information and facilitate exploitation within the optimization problem's search space. Further enhancing its adaptability, the SLLS incorporates a natural growth rule that generates varying learning efficiencies, adjusting with each iteration via the Sigmoid function. This adaptive feature balances the algorithm's exploration and exploitation capacities throughout the optimization process.

To validate the practical and theoretical merits of the SLLS, we subject it to rigorous testing using 60 widely recognized benchmark optimization problems and seven prominent engineering optimization challenges.

The primary contribution of this paper lies in the design of a novel metaheuristic algorithm, informed by the observation of snake locomotion and pheromone-based communication. By introducing a dynamic learning efficiency through the Sigmoid function, our algorithm adeptly balances exploration and exploitation, enhancing its problem-solving capabilities. The comprehensive evaluation on diverse benchmark problems and engineering challenges further substantiates the SLLS's value and versatility.

The subsequent sections of this paper are organized as follows: Section 2 provides an extensive review of the literature, categorically exploring the development and characteristics of various heuristic algorithms. Section 3 delves into the construction of the Snake Locomotion Learning Search algorithm (SLLS). Here, we comprehensively present the core concepts, primary operations, implementation procedure, control mechanisms, and the computational time complexity of the SLLS. Section 4 showcases the SLLS's efficacy and effectiveness through its application to 60 benchmark optimization problems and 7 renowned engineering optimization challenges. We conduct a comparative analysis of results obtained using various heuristic algorithms and also perform a sensitivity analysis related to the key parameters of the SLLS. Finally, Section 5 summarizes the principal conclusions drawn from our research and outlines potential avenues for future exploration. In Appendices A and B, readers will find a detailed presentation of the 30 benchmark optimization problems and the 7 well-established engineering optimization problems, respectively.

## 2. Literature review

In the following section, we will introduce a selection of heuristic algorithms in chronological order. Due to space constraints, we will omit discussion of most algorithm variants. These heuristics will be categorized into three distinct groups. The first group encompasses algorithms inspired by biological behaviors, excluding those derived from human social activities. The second group comprises algorithms rooted in human social activities. The third group encompasses heuristic algorithms inspired by physical phenomena and mathematical concepts.

Many algorithms in the first group were developed earlier in the past. Holland (1992)[3] proposed the Genetic Algorithm (GA). Kennedy and Eberhart (1995)[4] developed the Particle Swarm Optimization (PSO) algorithm. Storn and Price (1997)[5] introduced Differential Evolution (DE) for global optimization in continuous spaces. He et al. (2004)[6] advanced the PSO to tackle mechanical design optimization problems. Mezura-Montes and Coello (2005)[7] created a multimembered Evolution Strategy (ES) for global, nonlinear optimization problems. Karaboga (2005)[8] combined the foraging behavior of honey bees to form the Artificial Bee Colony (ABC) algorithm. Dorigo et al. (2006)[9] surveyed numerous applications of Ant Colony Optimization (ACO). Simon (2008)[10] presented Biogeography-Based Optimization (BBO). Yang and Deb (2009)[11] devised Cuckoo Search (CS) based on obligate brood parasitic behavior of certain cuckoos. Oftadeh et al. (2010)[12] proposed the Hunting Search (HS) meta-heuristic algorithm inspired by cooperative hunting habits of species. Yang (2010)[13, 14] developed the Firefly Algorithm (FA) and the Bat Algorithm (BA) for nonlinear design problems. Yang (2010)[15] compared the various heuristic algorithms proposed earlier. Through long-term development, these earlier heuristic algorithms, especially GA, PSO, DE, ES and ACO, are widely used today in tackling all kinds of optimization problems.

Over the past ten years, more heuristic algorithms in the first group have been proposed. Yang (2012)[16] presented the Flower Pollination Algorithm (FPA) which is inspired by the pollination process of flowers. Yang and Hossein (2012)[17] developed the Bat Algorithm (BA) to solve engineering optimization challenges. During this period, Mirjalili and his collaborators developed a series of heuristic algorithms including the Grey Wolf Optimizer (GWO) [18], the Moth-Flame Optimization (MFO) [19] algorithm, the Whale Optimization Algorithm (WOA) [20], the Salp Swarm Algorithm (SSA) [21], and the Multiobjective Salp Swarm Algorithm (MSSA) [21]. Azizyan et al. (2019)[22] established the Flying Squirrel Optimizer (FSO). Heidari et al. (2019)[23] created a population-based optimization paradigm called Harris Hawks Optimizer (HHO) based on the cooperative hunting procedures of Harris hawks. Abdollahzadeh et al. (2021)[24] developed the African Vultures Optimization Algorithm (AVOA), which mimics the foraging and navigational methods of African vultures. Recently, Zhao et al. (2022)[25] unveiled an artificial hummingbird algorithm (AHA), simulating the flight expertise and foraging techniques of hummingbirds.

In addition to the algorithms inspired by biological behaviors mentioned earlier, numerous studies have sought to replicate human behaviors and create efficient algorithms by simulating human behavior. Ray and Liew (2003)[26] proposed an optimization algorithm based on the information exchange occurring between and within societies. Lee and Geem (2005)[27] developed the Harmony Search (HS) algorithm for optimizing problems with continuous variables, using the concept of musical harmonies. Rao et al. (2011)[28] presented Teaching-Learning-Based Optimization (TLBO) for mechanical design optimization. He (2023)[29] proposed a novel swarm heuristic algorithm called the Medalist Learning Algorithm (MLA) which is inspired by the learning behavior of individuals in a group. He (2023)[29] and He and Cui (2023)[30] applied the MLA to truss structural optimization with natural frequency constraints and the configuration optimization of trusses. He and Cui(2023)[31] modified the original MLA to propose a new algorithm called the Multiscale Model Learning Algorithm (MMLA) and implement it to search the equilibrium state of a multi-tiered supply chain network. All of the aforementioned heuristic algorithms belong to the second group introduced earlier in this subsection. Because these algorithms typically draw their fundamental ideas from observations of human society, they are generally straightforward to understand and implement in practice.

In addition to the above two types of heuristic algorithms, many algorithms are derived from the observation and simulation of physical phenomena or the use of mathematical concepts. Kirkpatrick et al. (1983)[32] pioneered the Simulated Annealing algorithm (SA) where the idea of annealing in solid structures is applied in optimizing complex systems. Formato (2007)[33] introduced the Central Force Optimization algorithm (CFO) which takes inspiration from gravitational kinematics. Rashedi et al.(2009)[34] proposed the Gravitational Search Algorithm (GSA) based on the Newtonian law of gravity and mass interactions. Kaveh and Talatahari (2010)[35] developed the Charged System Search (CSS), which is based on the Coulomb law and the Newtonian laws of mechanics. Eskandar et al. (2012)[36] proposed the Water Cycle Algorithm (WCA) which is modeled on how streams and rivers run their course to the sea. Moghaddam et al. (2012)[37] presented the Curved Space Optimization (CSO), which makes use of the concept of space-time curvature. Hatamlou (2013)[38] developed the Black Hole Algorithm (BHA) and applied it to the clustering problem. Sadollah et al. (2013)[39] proposed the Mine Blast Algorithm (MBA) based on the concept of mine bomb explosion. Salimi (2015)[40] introduced the Stochastic Fractal Search algorithm (SFS), which is based on growth according to a mathematical fractal concept. Poonam and Vimal (2016)[41] presented the Passing Vehicle Search algorithm (PVS), which is based on the mathematics of vehicles passing on a two-lane highway. Kaveh and Dadras (2017)[42] developed the Thermal Exchange Optimization algorithm (TEO). Faramarzi et al. (2020)[43] presented the Equilibrium Optimizer algorithm (EO) inspired by the use of balance models to estimate both dynamic and equilibrium states. Zhao et al. (2020)[44] created the Artificial Ecosystem-based Optimization algorithm (AEO) based on the flow of energy in an ecosystem. Most of the above algorithms have been proved very effective to address some optimization problems.

Research into heuristic algorithms has been an active study area recently. Numerous new algorithms have been proposed and used in a variety of contexts. For example, Ma et al. (2021)[45] developed a multi-stage evolutionary algorithm for multi-objective optimization. Yildiz et al. (2022)[46] enhanced the classical grasshopper optimization algorithm by incorporating an elite opposition-based learning regime. Zhang (2023)[47] designed an Elite Archives-driven Particle Swarm Optimization (EAPSO) to improve the global searching ability of PSO. Ghasemi et al. (2023)[48] proposed a new geological inspired heuristic algorithm called Geyser inspired Algorithm (GEA). Rezaei et al. (2023)[49] presented the Geometric Mean Opimizer (GMO) for optimizing some well-known engineering problems. For a comprehensive review of existing heuristic algorithms, you can refer to Mohammed et al. (2023) [50], Yang (2023) [51], and Korani and Mouhoub (2021) [52].

## 3. Snake locomotion learning search (SLLS)

In this section, we will commence by elucidating the fundamental concepts underlying the SLLS. Following that, we will meticulously dissect the key operations, encompassing the implementation of serpentine and caterpillar movements, the establishment and updating of the list of visible spots, and the generation and utilization of learning efficiency. Subsequently, we will outline the implementation procedure of the SLLS. We will also delve into an analysis of the algorithm's computing time complexity and its control mechanism. To conclude this section, we will undertake a comparative examination of the SLLS in relation to other algorithms, focusing on their principal attributes.

### 3.1. Fundamental ideas

The fundamental concepts behind the Snake Locomotion Learning Search (SLLS) algorithm draw inspiration from the observation of snake motion within a region. These concepts rest on two pivotal aspects of snake behavior: two common types of snake locomotion and the use of pheromones to mark trails.

Serpentine movement, characterized by an "S" shape, is frequently employed by snakes for cruising in a region or swimming in water[53-55]. This mode of locomotion typically lacks a clear destination, serving the snake's goal of encountering an ideal spot for ambushing or resting.

Rectilinear movement (also called caterpillar movement), on the other hand, is the preferred choice when snakes need to approach prey or an ideal ambushing spot[53-55]. During rectilinear movement, snakes move slowly in a straight-line fashion, distinguished by vertical undulations.

These two modes of locomotion form the basis for designing the SLLS, enabling the creation of two types of search operations: random and deterministic.

Snakes, much like other animals, produce pheromones for various purposes, including attracting mates, defending territory, marking trails, and communication[56, 57]. In the SLLS, we assume that snakes mark potential spots for future exploitation. When a spot is marked, it becomes visible to nearby snakes, attracting them to explore and exploit the area around that spot.

Pheromones, however, evaporate over time. To account for this, we continuously update the list of visible spots. If a new spot exhibits a better fitness (or a smaller objective value, in the context of minimization problems) compared to some of the old visible spots, it replaces the least fit old spot, provided the list size exceeds a predefined limit.

Beyond the insights derived from snake behavior, we incorporate a natural growth rule into our algorithm's design. This entails generating a learning efficiency through the Sigmoid function. According to this rule[58, 59], learning efficiency starts at a low level and gradually increases during the initial stages of the learning process. As the process continues, the efficiency rises more rapidly. In the final phase, learning efficiency reaches a high level but increases at a much slower rate.

This characteristic of a typical learning process guides the SLLS. Learning efficiency controls the choice between serpentine movement and caterpillar movement throughout the solution process. As iterations progress, snakes are more likely to opt for caterpillar locomotion.

Additionally, the learning efficiency influences the amplitude of serpentine movement, causing it to decrease over time. This adjustment aims to focus the snake's energy on a more promising area as the algorithm evolves.

The connection between snake movement within a region and problem optimization is as follows: The snakes' search region corresponds to the feasible field (or search space) of an optimization problem. A feasible solution corresponds to a point within this space, referred to as a "spot" in this paper, reflecting snakes' searching behavior in an area. The attraction or fitness of a spot aligns with the objective function value of the associated solution.

As nearly all the optimization problems discussed in this paper involve minimizing an objective function, we assume that a spot's fitness is equivalent to the objective value in a minimization problem. In other words, a smaller solution equates to a higher fitness. Additionally, we assume that all objective values are greater than zero. If the minimal objective value is less than zero, we can introduce a sufficiently large positive constant to ensure that the resulting objective value is always greater than zero for all feasible solutions.

### 3.2. Key operations

#### 3.2.1. Serpentine locomotion

To simulate serpentine movement, we introduce key concepts. The distance from the starting point to the end spot of a serpentine movement is termed "amplitude" ($L_{\mathrm{A}}$). The end spot is the "foothold" or "final spot." As a snake moves, it forms an "S" or multiple "S" wave-like trail. Certain spots become pivotal for propulsion, known as "touch spots" (or points) with estimated fitness values.

To clarify these concepts, consider Figure 1. "1" is the starting point, while "2", "3", "4", and "5" represent touch points. The red curve in Figure 1 illustrates the serpentine movement from "1" to "5", with the straight-line distance between "1" and "5" being the amplitude.

Assuming "3" aligns precisely between "1" and "5", "8" (generating "2") becomes the midpoint between "1" and "3". Small squares labeled "6" and "7" are direction-guiding points.

The forward segment of the "S" bounded by "1" and "3" forms one half-circle, and the section between "3" and "5" creates the other half-circle.

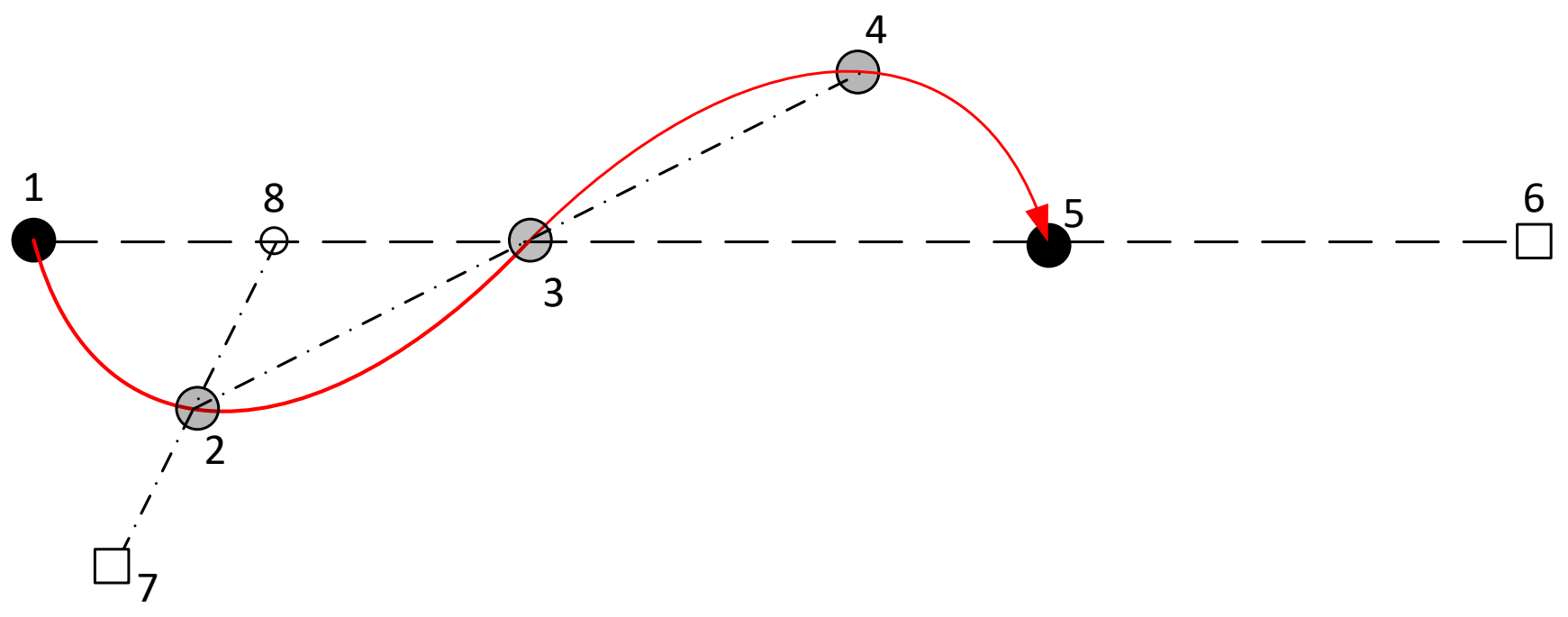

Fig. 1. Serpentine locomotion with four touch spots.

In later sections, we will delve into the analysis of serpentine locomotion, which may consist of more than two half-circles. For now, let's continue using Figure 1 to illustrate the stepwise procedure for generating a serpentine movement.

The stepwise procedure is as follows:

**Step 1**: Generate a random point within the feasible field, just as point "6" depicted in Figure 1.

**Step 2**: Determine the primary direction of the serpentine movement, which is the direction from the starting point "1" to the direction-guiding point "6", denoted as $DR_{1\to6}$.

**Step 3**: Along the primary direction, generate the foothold (the final touch point). The distance between the starting spot "1" and the foothold "5" is set to match the predefined amplitude of the serpentine movement. It's important to note that as the learning stage progresses in SLLS, the amplitude may gradually be reduced.

**Step 4**: If the generated foothold falls outside the feasible field, its individual components or elements should be adjusted. In other words, if a component exceeds its feasible range, its value should be substituted with the nearest boundary value. Notably, any adjustments made in this step may affect both the main direction and the amplitude of the actual serpentine movement. The adjusted result is retained.

**Step 5**: Utilizing the foothold and the starting spot, determine the middle touch point, analogous to point "3" in Figure 1. Likewise, pinpoint the middle point "8" between the starting spot "1" and the middle touch point "3".

**Step 6**: Similar to Step 1, generate a random point "7" within the feasible search space as a new direction-guiding point to ascertain the touch point "2."

**Step 7**: Along the direction $DR_{8\rightarrow7}$, ascertain the position of the touch point "2". Assume that the distance between touch point "2" and middle point "8" equals one-quarter of the actual amplitude determined in Step 4.

**Step 8**: Just as in Step 4, examine the feasibility of touch point "2" and make any necessary adjustments to ensure its feasibility.

**Step 9**: Using touch point "3" as the center, generate the mirror point "4" corresponding to touch point "2".

**Step 10**: Confirm the feasibility of mirror point "4" and, if needed, perform the necessary adjustments. The feasible point resulting from mirror point "4" becomes the final touch point "4".

**Step 11**: With the generated touch points "2", "3", "4", and "5", the serpentine movement, characterized by a random "S" shape, is fully formed. ■

The detailed mathematical operations required for the preceding procedure have been omitted here since they primarily involve basic arithmetic operations. Explicitly presenting these operations would make the procedure seem cumbersome and less accessible at first glance. We will provide a more generalized generation procedure for serpentine locomotion consisting of "$n$" half circles subsequently.

Serpentine locomotion may encompass more than just two half circles. For instance, Figures 2 and 3 illustrate two types of serpentine locomotion featuring 3 and 4 half circles, respectively. In the SLLS, each of these specialized forms of serpentine locomotion may be selected to realize the random search operation.

In Figure 2, points "2" to "7" represent the touch points. Points "8" and "10" are two auxiliary points randomly generated within the search space to facilitate the determination of points "7" and "2", respectively. Touch points "3" and "5" partition the line segment with "1" and "5" as endpoints into three equal segments. Point "9" serves as the midpoint between points "1" and "3". The touch point "2" positioned on the line connecting points "9" and "10", maintains the same distance from point "9" as point "1" does to point "9". Point "4" acts as the mirror point of point "2" with point "3" serving as the center. Point "6" is the mirror point of point "4", using point "5" as the reflection center.

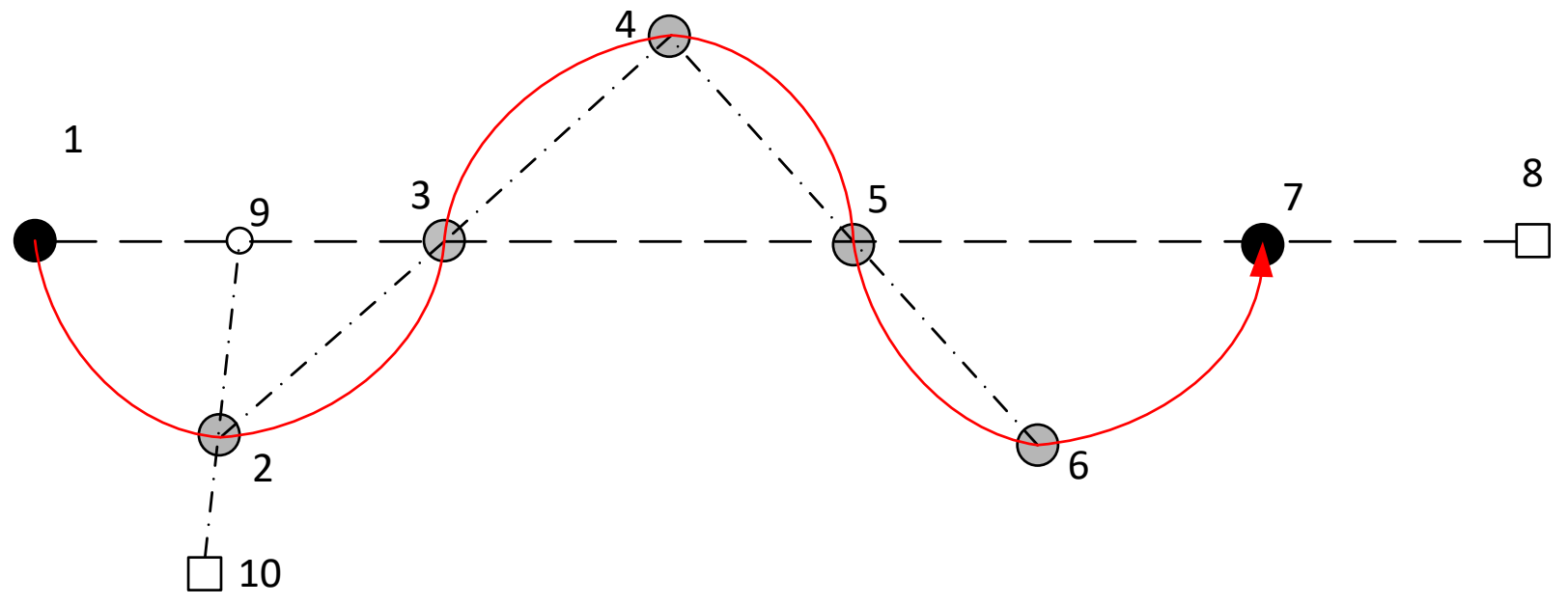

Fig. 2. Serpentine locomotion with 3 half circles and 6 touch points.

In Figure 3, points "2" to "9" represent the touch points. Points "3", "5", and "7" divide the line segment with points "1" and "9" as its two endpoints into four equal-length segments. Point "11" serves as the midpoint between points "1" and "3". Points "10" and "12" are two auxiliary points randomly generated within the search space to facilitate the determination of points "9" and "2", respectively. Point "2" lies on the line intersecting points "11" and "12." The distance from point "2" to point "11" is equivalent to the distance from point "1" to point "11". Points "2" and "4" form a pair of mirror points with point "3" as the center. Point "6" is the mirror point of points "4" and "8" with points "5" and "7" serving as the reflection centers, respectively.

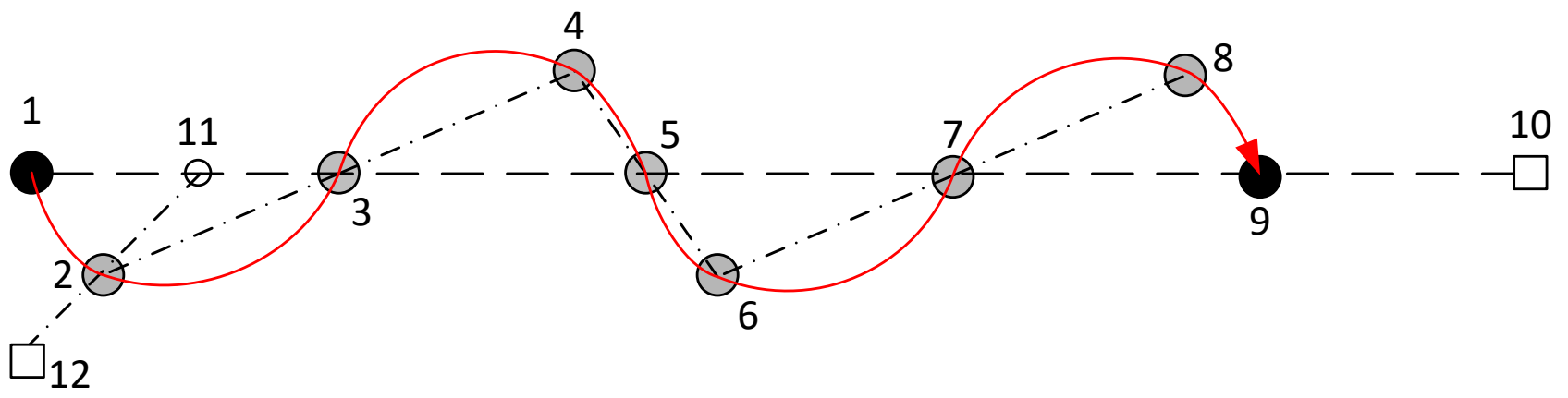

Fig. 3. Serpentine locomotion with 4 half circles and 8 touch points.

Before presenting the procedure for realizing serpentine locomotion with "$n$" half circles, where "$n$" is an integer equal to or greater than 1, we explain the related notations as follows: Let $X_i$ denote the coordinate vector of point "$i$" (or a spot in the $m$-dimensional search space), comprising a total of "$m$" elements, in the form $X_i = (x_{i,1}, x_{i,2}, \dots, x_{i,m})$. Here, "$m$" represents the number of decision variables. For simplicity in subsequent expressions, we will refer to $X_i$ as point "$i$". Let $[\underline{x}_k, \overline{x}_k]$ represent the range of feasible values for the element $x_{i,k}$, where $\forall i \in \mathrm{I}$. The set "I" is a set of integers that includes the indices of the points currently under consideration. In other words, we have $x_{i,k} \in [\underline{x}_k, \overline{x}_k]$, for all $i \in \mathrm{I}$ and $\mathrm{k} = 1,2, \dots, \mathrm{m}$.

Let $\Gamma(\cdot)$ represent the normalizing function, capable of transitioning a point located outside the feasible search space to a point within the feasible space. Suppose $X_j'$ is an arbitrary point, whether feasible or infeasible. We can apply the normalizing function to $X_j'$ to obtain its corresponding feasible point, denoted as $X_j$. The operation of the function $\Gamma(\cdot)$ is illustrated as follows:

$$X_j = \Gamma(X_j'), \tag{1.1}$$

such that

$$x_{j,k} = \begin{cases} \underline{x}_k & , if\ x_{j,k}' < \underline{x}_k; \\ x_{j,k}' & if\ \underline{x}_k \le x_{j,k}' \le \overline{x}_k; \forall \mathrm{k} \\ \overline{x}_k & , if\ x_{j,k}' > \overline{x}_k. \end{cases} \tag{1.2}$$

Let $l_{i,j}$ be the distance between spots $i$ and $j$ that equals $\|X_j - X_i\| = \sqrt[2]{\sum_k (x_{j,k} - x_{i,k})^2}$. If only the box constraints on variable values are considered, the above eq. (1) is equivalent to solving the following minimal problem:

$$X_j = \underset{X_k \in X}{\operatorname{argmin}} \|X_k - X_j'\|, \tag{2}$$

where the feasible region, X, is defined by the box constraints, which consist of interval constraints on the components of a variable.

For the sake of convenience in subsequent analysis, we provide the following definitions: As mentioned earlier, $DR_{i \to j}$ represents the direction from point "$i$" to point "$j$". The vector resulting from $\frac{1}{2}(X_i + X_j)$ is the coordinate vector of the middle point between node "$i$" and node "$j$". Similarly, the vector obtained by calculating $X_i + r(X_j - X_i)$ is the coordinate vector of a point on the line connecting nodes "$i$" and "$j$", with the distance between this point and "$i$" denoted as $r \times l_{i,j}$, where $r \in [0,1]$. Likewise, $X_i + \frac{L}{l_{i,j}}(X_j - X_i)$ is computed to determine the coordinate vector of a point at a distance of "$L$" from node "$i$" along the direction $DR_{i\to j}$. Finally, by calculating $2X_j - X_i$, we derive the coordinate vector of the mirror point of node "$i$" with node "$j$" as the center point.

Let $F(X_i)$ or simply $F_i$ represent the fitness value or objective function value of $X_i$. As previously mentioned, we assume that $F(X_i) \geq 0$, and the optimization problem involves minimizing $F(X_i)$, with or without constraints.

Now, let's explore how to implement serpentine locomotion consisting of "$n$" half circles.

Suppose the current location of the snake, which is also the starting point, is denoted as $X_1$. Here, "$n$" represents the number of half circles. Building upon the previous analysis involving Figures 1 to 3, it's evident that there are $2n$ touch spots associated with serpentine locomotion featuring "$n$" half circles. These touch spots collectively form a list: $\{X_2, X_3, X_4, \ldots, X_{2n+1}\}$. These spots are positioned along the serpentine trail in the order presented in the list. We can further categorize these touch spots into three sub-groups: $\{X_{2n+1}\}$, $\{X_3, X_5, X_7, \ldots, X_{2n-1}\}$ and $\{X_2, X_4, X_6, \ldots, X_{2n}\}$. The single spot in the first sub-group corresponds to the end spot of the serpentine locomotion. The spots in the second group lie on the line connecting the starting spot $X_1$ and the end spot $X_{2n+1}$, effectively dividing the line between $X_1$ and $X_{2n+1}$ into "$n$" segments of equal length. The spots in the third group are touch points positioned outside the line between $X_1$ and $X_{2n+1}$, contributing to the formation of multiple "S" shapes.

Let's outline the stepwise procedure for realizing serpentine movement with "$n$" half circles:

**Step 1:** Generate the end spot $X_{2n+1}$ (referred to as the foothold or the final spot of the serpentine movement). First, randomly generate a point $X_{2(n+1)}$ in the feasible search space as an auxiliary point to determine the foothold. Then, calculate $X_{2n+1}$ as follows:

$$X_{2n+1} = \Gamma\left[X_1 + \frac{L_A}{l_{1,2(n+1)}}(X_{2(n+1)} - X_1)\right]. \tag{3}$$

In equation (3), $L_A$ represents the given amplitude of serpentine locomotion. After applying the adjustment operation $\Gamma(\cdot)$, the actual amplitude should be updated as follows:

$$\tilde{L}_A = \|X_{2n+1} - X_1\|. \tag{4}$$

**Step 2:** Determine the spots in group $\{X_3, X_5, X_7, \ldots, X_{2n-1}\}$. For a spot $X_i$ in this group, it can be calculated using the following formula:

$$X_i = X_1 + \frac{i-1}{2n}(X_{2n+1} - X_1). \tag{5}$$

**Step 3:** Generate a random point $X_{2(n+2)}$ in the feasible search space as an auxiliary point to determine $X_2$. Then, obtain the middle point $X_{2n+3}$ between points $X_1$ and $X_3$ as follows: $X_{2n+3} = \frac{1}{2}(X_1 + X_3)$. Finally, determine $X_2$ using the equation:

$$X_2 = \Gamma\left[X_{2n+3} + \frac{l_{1,2n+3}}{l_{2n+3,2(n+2)}}(X_{2(n+2)} - X_{2n+3})\right]. \tag{6}$$

**Step 4:** Determine the points in group $\{X_4, X_6, \ldots, X_{2n}\}$. For a spot $X_i$ in this group, it can be calculated as follows:

$$X_i = \Gamma[2X_{i-1} - X_{i-2})]. \tag{7}$$

Note that the points in $\{X_4, X_6, \ldots, X_{2n}\}$ should be calculated in the order as given in the list, starting with $X_4$ and progressing until the last point $X_{2n}$.

**Step 5:** By placing all the obtained points in a list $\{X_2, X_3, X_4, \ldots, X_{2n+1}\}$, we arrive at the final set of touch points that constitute the serpentine movement. ■

### 3.2.2. Caterpillar locomotion

Caterpillar locomotion is a common observation in snakes. When approaching a potential spot where prey may be found, a snake moves slowly and quietly along a straight-line path, resembling the undulating movement of a caterpillar on a twig. To simulate this specific movement within an optimization process, we will use the caterpillar movement depicted in Figure 4 as an illustrative example to demonstrate how to realize a caterpillar movement with four designated touch points. These touch points represent the locations that the snake will assess, meaning that the fitness values of these points will be estimated. Later, we will explain the more generalized form of caterpillar movement involving a predetermined number of touch points.

In Figure 4, as a snake approaches the target spot, indicated by a red triangle, it only checks four specific points—points 2 to 5—to estimate their fitness values, even though there may be multiple touch points along the way. This selection of four check points is strategically unevenly dispersed along the straight trail, as the snake becomes more cautious and meticulous when nearing its target spot.

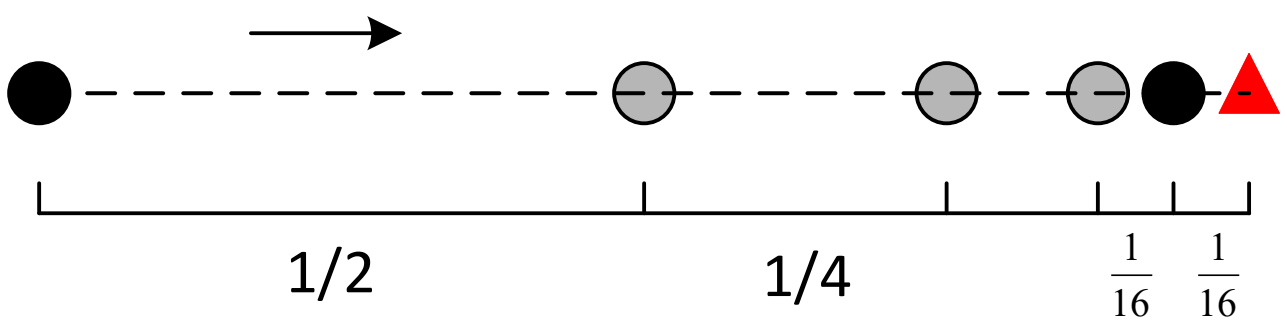

Fig. 4. Caterpillar locomotion with four touch points.

In Figure 4, spot 1 serves as the starting point for the snake's movement, while spot 6, marked with a red triangle, represents the target location. Spots 2, 3, 4, and 5 are designated as the four touch points used in this movement. The distance between the starting spot and the target spot is referred to as $L_{CL}$, known as the rectilinear distance of caterpillar locomotion.

In Figure 4, we introduce a location parameter $r_{CL}$ with a value of 0.5, referred to as the rate for demarcation. The first touch point, point 2, is located at a distance of $r_{CL}L_{CL}$ from the starting point 1. In this context, point 2 coincides with the midpoint between points 1 and 6. The second touch point, point 3, is situated at $r_{CL}L_{CL} + r_{CL}(1 - r_{CL})L_{CL}$ from the starting point 1. Similarly, point 3 corresponds to the midpoint between touch point 2 and the target spot 6. The third touch point, point 4, is the midpoint between points 3 and 6 in Figure 4. The final touch point, point 5, evenly divides the path between points 4 and 6 in Figure 4. This last touch point serves as the foothold for the caterpillar movement.

It's important to note that the choice of touch points and their arrangement along the path is not unique. On one hand, readers can explore various methods for positioning touch points along the approach path using different demarcation rates ($r_{CL}$). However, the value of $r_{CL}$ should remain within the interval (0, 1). On the other hand, the number of touch points checked during a caterpillar movement can also lead to variations in the movement, even when $r_{CL}$ remains the same.

Assuming the feasible field is convex, the touch points generated as described above will yield feasible solutions because both the starting point and the target point are located within a convex feasible field.

Upon completing a caterpillar movement, the snake selects the last touch point, such as point 5 shown in Figure 4, as its final foothold. This choice is made to avoid redundant computations because the target spot is typically a solution that has already been checked and stored in the list of visible spots.

Now, let's expand upon the special caterpillar movement with four touch points to caterpillar movements with $n$ touch points, while also considering a specified rate for demarcation, $r_{CL}$, which is restricted to the interval (0.0, 1.0).

Suppose we have $X_1$ as the starting spot, and $X_V$ is the visible spot chosen from the list of visible spots to serve as the target point. Our goal is to determine a series of touch points $\{X_2, X_3, \dots, X_n\}$ that will be checked during a caterpillar locomotion. Let $X_i$ represent a typical touch point to be checked from the set $\{X_2, X_3, \dots, X_n\}$. We can use the following sequence of equations to compute the touch points one by one, following the order presented in the list $\{X_2, X_3, \dots, X_n\}$:

$$\begin{cases} X_2 = X_1 + r_{CL}(X_V - X_1); \\ X_3 = X_2 + r_{CL}(X_V - X_2); \\ \quad\vdots \\ X_i = X_{i-1} + r_{CL}(X_V - X_{i-1}); \\ \quad\vdots \\ X_n = X_{n-1} + r_{CL}(X_V - X_{n-1}). \end{cases} \tag{7}$$

By a series of simple algebraic manipulations, we have the following equivalent uniform equation to calculate $X_i$ directly:

$$X_i = X_1 + \left[\sum_{k=1}^{i-1} (-1)^{k-1} C_{i-1}^{k}\, {r_{CL}}^{k}\right] (X_V - X_1), \forall i = 2,3,\dots,n. \tag{8}$$

In eq. (8), $C_n^k$ stands for the combination $C_n^k = \frac{n!}{k!(n-k)!}$ and k! is the factorial defined as $k! = 1 \times 2 \times 3 \times \dots \times k$.

The choice of different combinations of $r_{CL}$ and $n$ for a caterpillar locomotion can significantly affect the performance of the SLLS. We will delve into the impact of these choices in the numerical analysis section. Unless explicitly specified otherwise, we will use a value of 0.5 for $r_{CL}$ in most cases throughout this paper.

### 3.2.3. Updating the list of visible spots

Some snakes release pheromones to mark trails and communicate with others. To simulate this mechanism in the SLLS and aid the optimization process, we will establish a list of visible spots. A visible spot refers to a solution in the optimization problem that has been assessed by a snake during its movement.

As previously discussed in subsections 3.2.1 and 3.2.2, all the touch points within the serpentine and caterpillar movements are verified solutions. However, considering all of them as visible points would result in an excessively large quantity. Based on this observation, we will selectively designate certain touch points with high fitness values as visible points and arrange them in an ordered list.

We will denote the list of visible spots as K. We assume that the size of K is predetermined and denoted by "$m$". The reason for predefining a size for the list of visible spots is that, while every touch spot may be marked with a certain amount of pheromone, we assume that only some of the spots with better fitness can effectively attract the attention of snakes for further exploration.

Let the list be represented as $\mathrm{K}=\{X_1^{\wedge}, X_2^{\wedge}, \ldots, X_m^{\wedge}\}$, where $X_k^{\wedge}$ for $k=1,2,..,\mathrm{m}$ represents a typical element of this list. To form an ordered list, the elements in K need to meet the following condition:

$$F(X_k^{\wedge}) \leq F(X_{k+1}^{\wedge}), \forall k=1,2,m-1. \tag{9}$$

The list of visible spots involves two primary operations: adding a new spot to the list and randomly selecting a visible spot from the list.

Firstly, let's discuss the process of adding a spot to the list, denoted as K. According to the earlier definition, the fitness value of any visible point in the list is greater than that of all the other touch points that have already been checked but are not in the list. The list of visible spots should be updated whenever a new touch point is checked. In terms of pheromone representation, visibility reflects the intensity of pheromone. For simplicity in subsequent analysis, we assume that the quantity of pheromone at a spot is proportional to the fitness of that spot. To facilitate the simulation, we also assume that the pheromone of a currently used visible point remains constant and does not evaporate throughout the entire search process.

To add a spot to the visible spot list, we must compare its fitness value with those associated with the spots already in the list. We represent this addition operation using an abstract function as follows:

$$\mathrm{K} := \Omega(\mathrm{K}, X_i). \tag{10}$$

In the abstract function described in Equation 10, $X_i$ represents the spot that we intend to add to the list K. The operational logic of the function $\Omega(K, X_i)$ is as follows: If $F(X_i) < F(X_1^{\wedge})$, then $X_1^{\wedge} := X_i$ and $X_k^{\wedge} := X_{k+1}^{\wedge}$ for all $\mathrm{k}=2,3,\ldots,m-1$; else if $F(X_k^{\wedge}) \leq X_i < F(X_{k+1}^{\wedge})$ for $k=2,3,m-1$, then $X_h^{\wedge} := X_h^{\wedge}$ for $h=1,2,k$, $X_{k+1}^{\wedge} := X_i$ and $X_j^{\wedge} := X_{j+1}^{\wedge}$ for all $j=\mathrm{k}+2,\mathrm{k}+3,\ldots,m-1$; else, $F(X_i) \geq F(X_m^{\wedge})$, then $X_k^{\wedge} := X_k^{\wedge}$ for all $\mathrm{k}=1,2,3,\ldots,m$.

During the initial phase of implementing SLLS, the actual size of the list K may be smaller than the predetermined size, denoted as $m$. In such cases, a new spot is directly added to the list in a suitable position to ensure that condition (9) is satisfied, and no existing visible spot is removed from the list. Since this procedure closely resembles the previously described operation logic, we will not repeat the explanation here.

In addition to adding a spot to the visible list, another operation associated with the list of visible spots is the selection of a visible spot from the list as the target point for a snake to execute a caterpillar movement. In this context, the roulette wheel selection method is employed. If the size of the visible list is sufficiently large, tournament selection may be considered. The process of selecting a visible point from the visible set is outlined below. As previously mentioned, we assume that all fitness values are non-negative, and we are currently dealing with a minimization problem. In a minimization problem, the smaller the solution, the higher the fitness value associated with the corresponding spot. Initially, we calculate the sum of the inverses of all the fitness values as follows:

$$\Sigma\mathrm{F} = \sum_{k=1}^{m} \frac{1}{F(X_k^{\wedge})}. \tag{11}$$

The probability of selecting the spot $X_k^{\wedge}$ from K is obtained as follows:

$$\Pr(X_k^{\wedge}) = \frac{(1/F(X_k^{\wedge}))}{\Sigma\mathrm{F}}, \text{for } k=1,2,\ldots,m. \tag{12}$$

As roulette wheel selection is employed here, we generate a random number, denoted as $Ran$, from a uniform distribution within the interval [0, 1]. Define $\Sigma\mathrm{P}_0 = 0$ and compute the accumulated probability as follows:

$$\Sigma \mathrm{P}_j = \sum_{k=1}^{j} Pr(X_k^{\wedge}), j = 1,2,\dots,m. \tag{13}$$

The rule to select a visible spot is as follows: If $\Sigma \mathrm{P}_{j-1} < Ran \leq \Sigma \mathrm{P}_j, j = 1,2,\dots,m$, then $X_j^{\wedge}$ will be chosen as the targeting spot.

**3.2.4. Generation and utilization of the learning efficiency**

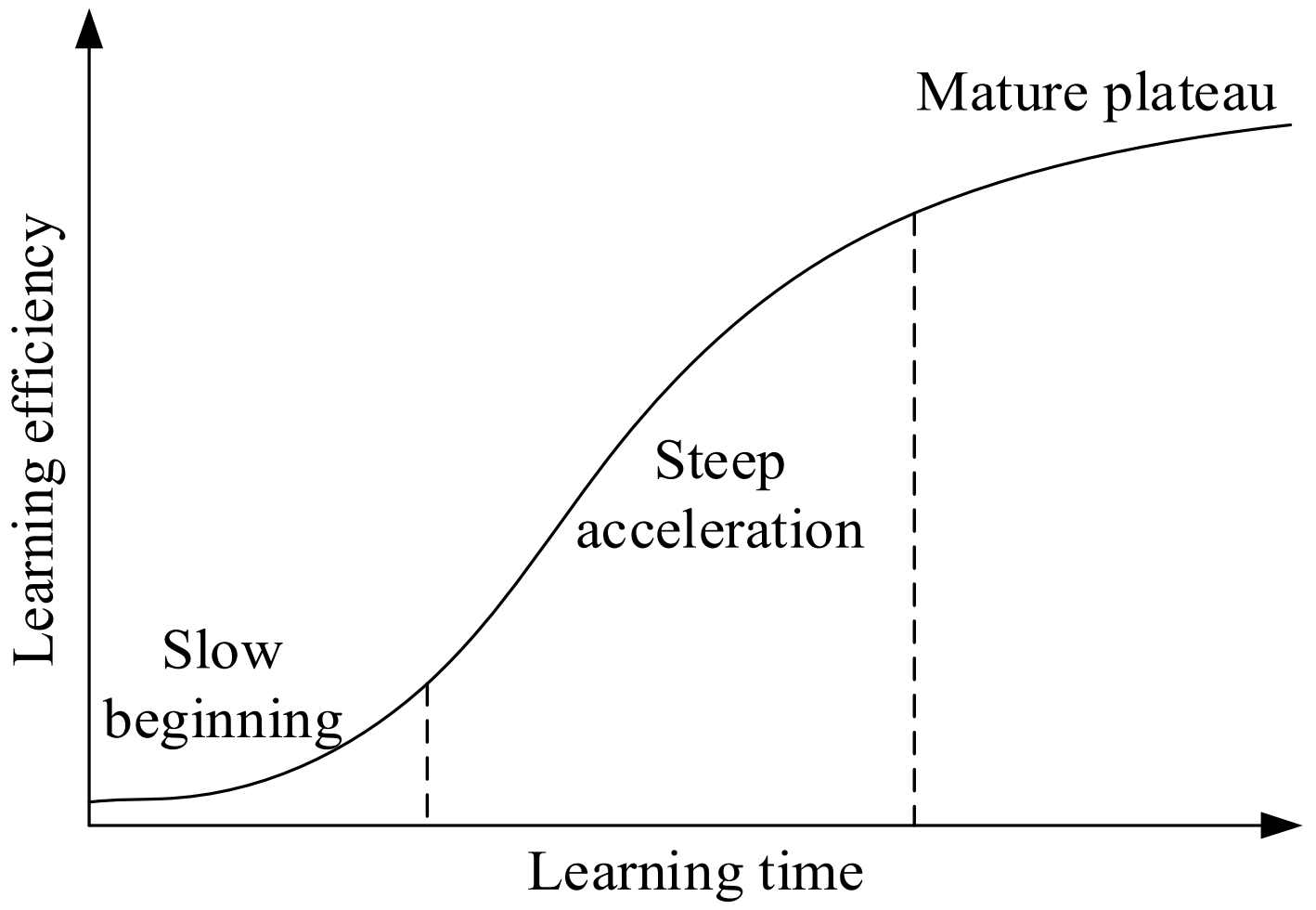

Fig. 5. Natural learning curve with three phases.

In the preceding three subsections, we introduced concepts inspired by snake behavior. In this subsection, we will introduce a concept drawn from biological growth.

For animals and humans, the process of learning typically follows a natural growth curve, often referred to as the Logistic growth curve, as illustrated in Figure 5 [58, 59]. Learning efficiency tends to change over time during the learning period. Initially, learning efficiency is low and increases gradually. As the learning process advances into its middle phase, efficiency experiences rapid improvement. In the final phase of learning, efficiency remains high but increases very slowly.

The above pattern of learning efficiency throughout a learning process can serve as a foundation for balancing the exploitation and exploration capabilities of a heuristic algorithm, such as MLA and MMLA [29-31]. Generally, in the early stages of searching for global optima using a heuristic algorithm, the emphasis should be on exploring the entire search space. However, as the algorithm accumulates experience and gains more information from past efforts, it should utilize the known global information to exploit the more promising regions of the search space, given its limited search resources.

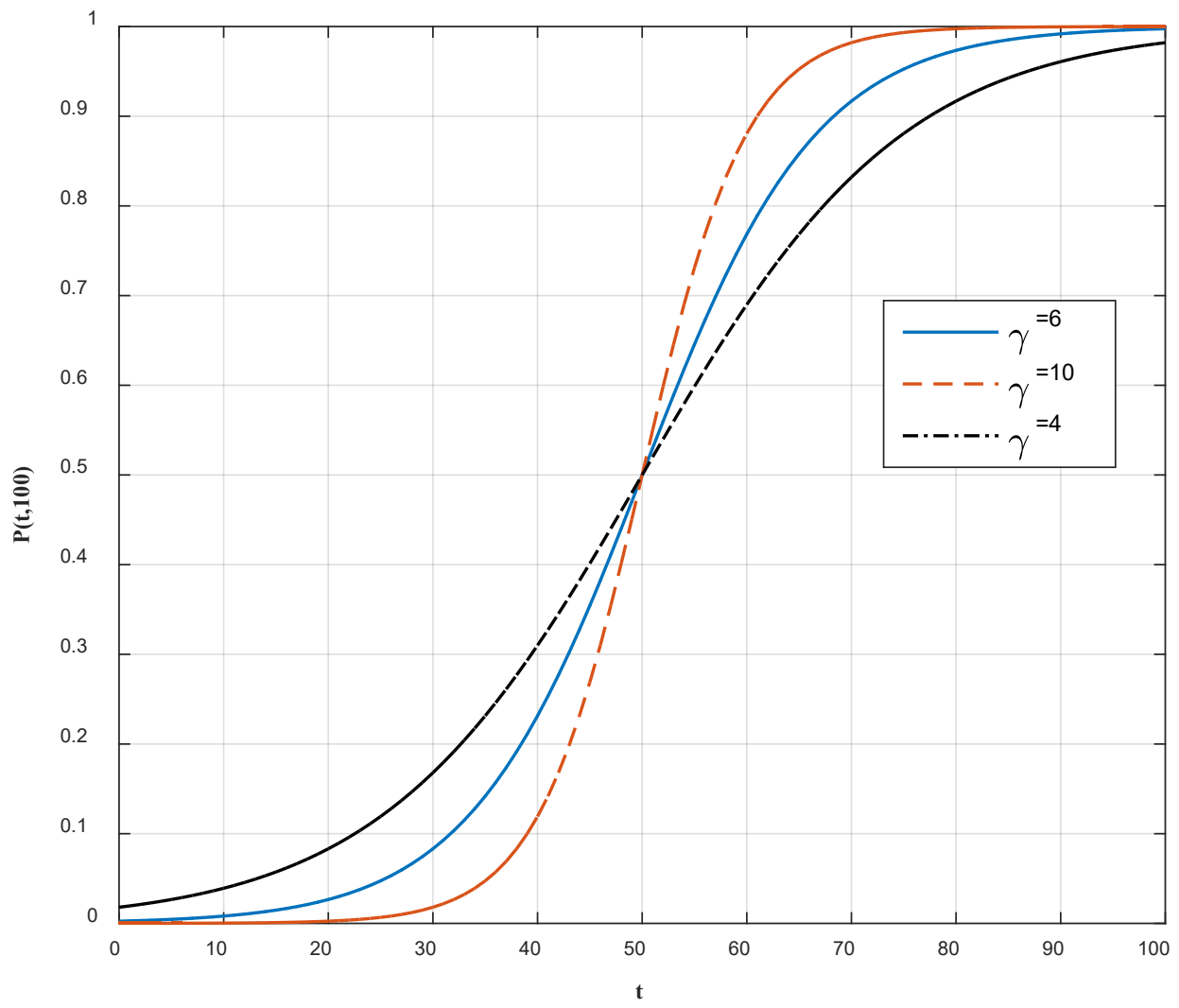

Fig. 6. Learning efficiency over the learning period.

In this paper, we will generate the learning efficiency using the Sigmoid function.

Let's assume that the process of executing a serpentine or caterpillar movement represents one movement for a snake, which corresponds to a single iteration of the SLLS. As described earlier regarding these two modes of snake locomotion, during one such movement, a snake will examine a series of touch points and ultimately select the last touch point as its foothold at the conclusion of the movement. For the purposes of our discussion, let's assume that during one movement (or one iteration of the SLLS), the learning efficiency remains constant for all the snakes. All snakes simultaneously commence and conclude a movement tailored to their specific needs.

Now, let's delve into the process of generating the learning efficiency. Assuming that the total number of movements for a snake (or the total number of iterations of the SLLS) is denoted as $T$, and the current iteration is represented by $t$, we can calculate the learning efficiency, denoted as $P(t)$, using the following equation:

$$P(t) = \frac{1}{1 + e^{\frac{2\gamma}{T}\left(\frac{T}{2} - t\right)}} \tag{14}$$

where $\gamma$ is a given positive constant referred to as the radius of the truncated range of the natural growth function, i.e. the Sigmoid function. As $t$ increases, $P(t)$ will start at a low value, gradually rise, then enter a rapidly increasing phase, remain at a steady plateau at a high level for a while, and finally converge to 1. Figure 6 depicts the variation of $P(t)$ with respect to $t$ when $T = 100$ and $\gamma$ equals 4, 6 and 10, respectively. The curves in Figure 6 show that for a given $T$, the bigger the value of $\gamma$ is, the steeper the rising phase will be.

The concept of learning efficiency plays a dual role in the SLLS. Firstly, it serves to balance the exploration and exploitation capabilities of the SLLS. Secondly, it is utilized to adapt the amplitude of serpentine locomotion as the number of SLLS iterations increases.

As previously mentioned, the exploration and exploitation capabilities of the SLLS are primarily represented by the serpentine and caterpillar movements, respectively. Achieving a balance between these two capabilities involves selecting the appropriate movements for a snake during the optimization process. As the number of iterations increases, the corresponding learning efficiency is computed using equation (14). In this paper, we assume that for a given iteration $t$, the probability of a snake opting for a serpentine movement is $1 - P(t)$, while the probability of conducting a caterpillar movement is $P(t)$. When deciding on a specific movement for a snake, a random number is generated from a uniform distribution spanning the interval [0,1]. If this random number is less than $P(t)$, the snake will perform a caterpillar movement; otherwise, it will engage in a serpentine movement.

Another application of learning efficiency is to adapt the expected amplitude of serpentine movement. Let's assume the current iteration is represented by $t$. The initial amplitude of serpentine movement is denoted as $L_A(0)$, while the amplitude associated with iteration $t$ is represented as $L_A(t)$. Additionally, there exists a predetermined minimal amplitude for serpentine movement denoted as $L_A^{min}$. The calculation of $L_A(t)$ is as follows:

$$L_A(t) = L_A(0) - \left[L_A(0) - L_A^{min}\right]P(t), \forall t \in \{1,2,\dots,T\}. \quad (15)$$

As evident, with the progression of iteration $t$, the amplitude of serpentine locomotion gradually decreases due to the increasing learning efficiency $P(t)$ throughout the SLLS iterations. It's worth noting that by diminishing the amplitude, SLLS directs the snakes' search efforts towards more promising areas.

### 3.3. Implementation procedure of SLLS

With the concepts and operations explained in subsections 3.1 and 3.2, the step-wise implementation procedure of the Snake Locomotion Learning Search (SLLS) algorithm is presented as follows:

**Step 1: Initialization**

Initialize various parameters, including the number of snakes $n_S$, the maximal number of iterations $T$, the initial amplitude of serpentine locomotion $L_A(0)$, the truncated range of the natural growth function $\gamma$, the number of half circles in a serpentine movement, the minimal amplitude of serpentine movement $L_A^{min}$, the number of touch points in a caterpillar movement, and the predetermined size of the list of visible points $m$. Set the current iteration $t$ to 1.

**Step 2: Generate Snakes and Initial List of Visible Points**

Generate a swarm of snakes with the given size $n_S$ and initialize the initial list of visible points. To generate a snake, create a random spot in the search space with a unique serial number. The initial list of visible spots consists of the current positions of snakes with relatively high fitness values compared to those outside the list.

**Step 3: Update Learning Efficiency and Amplitude**

Update the learning efficiency $P(t)$ and amplitude $L_A(t)$ for the current iteration $t$ using equations (14) and (15), respectively.

**Step 4: Conduct Snake Movements**

Conduct movements for each snake associated with iteration $t$. Perform the following operations one by one for all the snakes in the swarm:

**Step 4.1: Choose Movement Mode**

Generate a random number for each snake from a uniform distribution in the interval [0, 1]. If the random number is greater than $P(t)$, proceed to Step 4.2; otherwise, move to Step 4.3.

**Step 4.2: Serpentine Movement**

Execute the serpentine movement for the snake following the explanation in subsection 3.2.1. Estimate the fitness values for all the touch points encountered in this movement and consider updating the list of visible spots with these touch spots.

**Step 4.3: Caterpillar Movement**

Perform a caterpillar movement for the current snake following the explanation in subsection 3.2.2. Estimate the fitness values for all the touch points encountered in this movement and consider adding these touch points to the list of visible spots.

**Step 5: Termination Criteria**

Check whether to terminate the algorithm. If $t$ is equal to $T$ or the absolute difference between the best and worst fitness values associated with the visible spots is less than a given positive constant $\delta_F$, terminate the algorithm and output the visible spot with the best fitness value. Otherwise, increment $t$ by 1 and return to Step 3.■

**3.4. Control mechanism and time complexity**

Figure 7 presents the control flow chart of the Snake Locomotion Learning Search (SLLS) algorithm, corresponding to the implementation process described in Subsection 3.3. The exploration and exploitation capacities of the SLLS are primarily manifested through serpentine locomotion and caterpillar locomotion of the snakes, respectively. Unlike some other heuristic algorithms, such as Ant Colony Optimization (ACO), which combine these two capacities into a single operation, the SLLS separates them into distinct operations that can only be carried out one at a time. The choice of locomotion is determined randomly based on a learning efficiency specific to the current learning stage (or, in other words, the current iteration of the SLLS). Historical experience is accumulated using an ordered list of visible spots, which are touch points with higher fitness values than all other checked points outside the list. This list is updated whenever a new spot is checked. The SLLS will terminate if the maximum number of iterations is reached or if the difference between the highest and lowest fitness values of the visible spots is smaller than a predefined small positive constant. During the implementation process of the SLLS, each of the snakes will undergo a randomly chosen locomotion in synchrony.

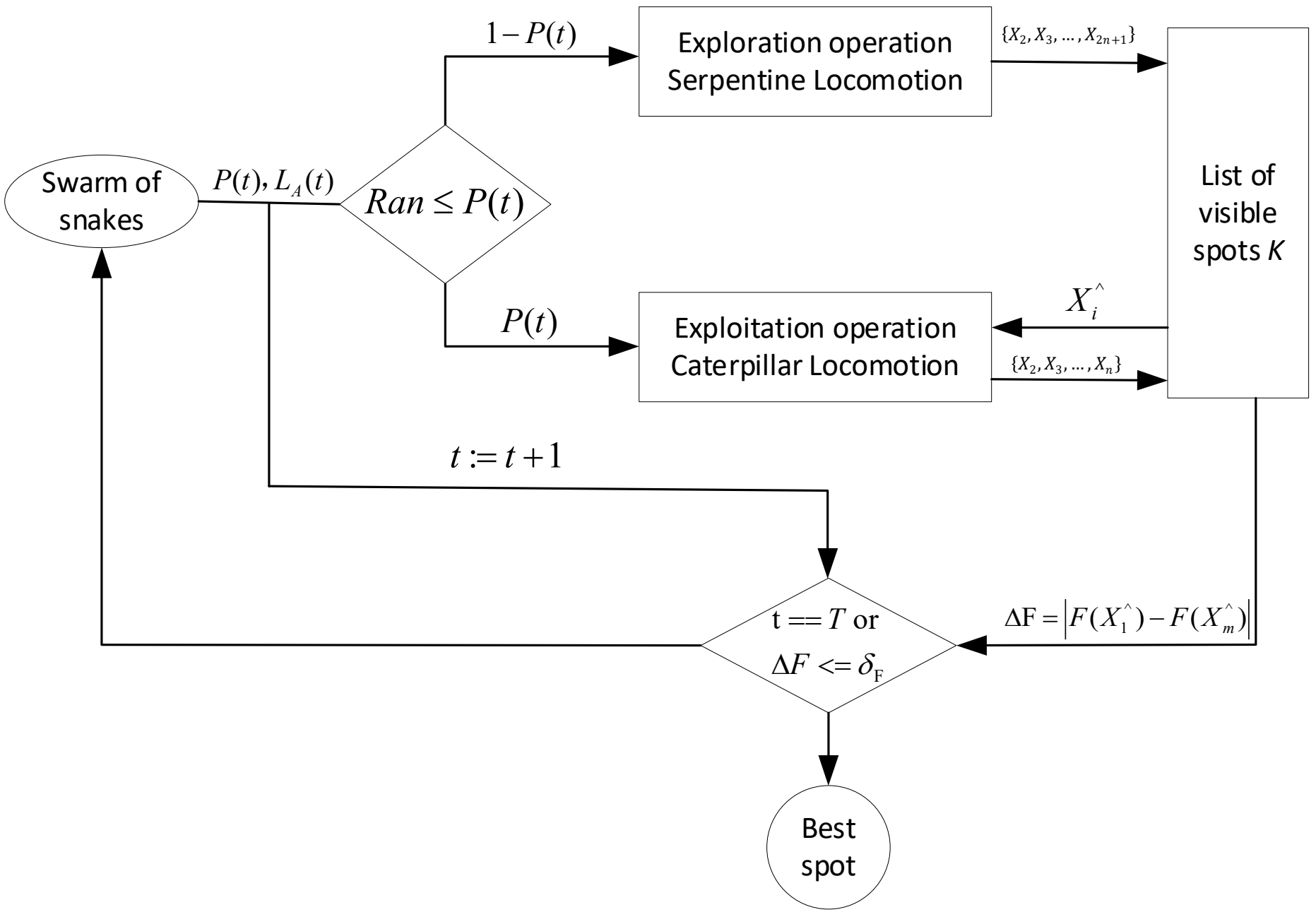

Fig. 7. Control flow chart of the SLLS.

With the above brief introduction to the control flow chart of the SLLS, we can roughly deduce the total number of function estimations (NFE). Let $n_{CM}$ and $n_{HC}$ denote the number of touch points in a caterpillar movement and the number of half circles in a serpentine movement, respectively. Clearly, the number of touch points in a serpentine movement will be $2n_{HC}$. Considering the symmetry of the curve of learning efficiency about the point with the coordinates $(T/2, 0.5)$, the chance or the number of executing each locomotion will be statistically equal for a snake. If $n_{CM}$ and $2n_{HC}$ are not equal, we have:

$$\text{NFE} \cong n_S \times T \times (n_{CM} + 2n_{HC})/2. \tag{16}$$

However, if $n_{CM}$ and $2n_{HC}$ are equal, the value of NFE can be obtained by:

$$\text{NFE} = n_S \times T \times n_{CM}. \tag{17}$$

Please note that in equations (16) and (17), $n_S$ and $T$ represent the number of snakes and the maximum number of iterations, respectively, as defined earlier.

Sometimes, we only need to compare the Number of Total Movements (NTM) conducted by all snakes when applying SLLS to an optimization problem many times. In this case, NTM can be calculated as:

$$\text{NTM} = n_S \times T. \tag{18}$$

Regarding the calculations associated with the determination of the touch points for a movement, they involve simple algebraic operations. Therefore, the time complexity of the SLLS can be reasonably estimated using the indices provided in equations (16-18). It's important to note that this estimation corresponds to the worst case, as the algorithm may terminate prematurely before reaching the maximum number of iterations.

### 3.5. The comparison with other heuristic algorithms

Nowadays, many new heuristic algorithms are continuously proposed to address various optimization problems, particularly those of considerable complexity in real-world scenarios. Evaluating the efficiency, effectiveness, and uniqueness of these algorithms can be a challenging task within this domain.

As per the "No Free Lunch Theorem"[2], no single algorithm can universally outperform all others across every conceivable problem scenario. In essence, if an algorithm excels in certain scenarios, it may exhibit less favorable performance in others. Therefore, when a heuristic algorithm demonstrates efficient and effective solutions for specific problem types, its practical and theoretical value should be acknowledged. In the following section, we will apply the SLLS to 67 benchmark problems to substantiate its value in accordance with the aforementioned observation.

Considering the uniqueness of an algorithm often involves comparing its optimization mechanisms with those of other algorithms, especially established heuristic ones. Table 1 provides a comparison between the SLLS and five classical heuristic algorithms, namely the Medalist Learning Algorithm (MLA), Genetic Algorithm (GA), Particle Swarm Optimization (PSA), Ant Colony Optimization (ACO), and Simulated Annealing (SA). This comparison spans six dimensions that capture the primary characteristics of these algorithms. The analysis reveals that the SLLS distinguishes itself not only in its fundamental operations for achieving exploration and exploitation capabilities but also shares some common features with existing algorithms.

The SLLS's uniqueness stems from its adoption of two distinct modes of snake locomotion. It shares the common trait of generating learning efficiency through a natural learning rule with the MLA. The use of pheromone in SLLS's construction is akin to ACO. The concept of a list of visible spots resembles the medalists in MLA. Furthermore, if the SLLS employs only one visible spot, this adjustment aligns with the use of a global best particle in PSO. It is conceivable to extend this comparison to a broader range of algorithms, but due to space limitations, we refrain from conducting a more comprehensive analysis here.

**Table 1**

Comparison of main mechanisms of several heuristic algorithms.

| Name | Inspired by | Exploration | Exploitation | Evolvement in time | Termination criteria | Balance mechanics |
|---|---|---|---|---|---|---|
| SLLS | Snake's locomotion | Serpentine locomotion | Caterpillar locomotion making use of pheromone information | Different stage (or iteration) is associated with a specialized movement | Maximal iterations or fitness values of all visible spots are near to each other | Learning efficiency based on natural growth rule |
| MLA | Individual's learning behavior in a group | Randomly changing a learning features in a neighborhood | Learning from the medalists or owns historical experience | Update the medalists in different learning stage | Maximal learning stage or medalists with the same performance | Learning efficiency based on natural growth rule |
| GA | Darwinian principle of survival of the fittest | Mutation operations and crossover operation | Crossover operation and selection operation | New generation obtained by implementing selection, crossover, and mutation | Maximal number of generation or no improvement on the best individual for a series generation | Using the proper parameters such as mutation and crossover rates |
| PSO | Birds' flying behavior in a flock | Move in the inertia direction | Move to the global and the local best location | Change the speeds and directions of movements of particles with time | Maximal number of iterations or particles flock together in a very narrow space | Control by proper parameters |
| ACO | Ants' foraging behavior | According the density of pheromone to choose a branch randomly | More likely go along the path with best visibility and high density of pheromone | Mark the trail with pheromone the quantity of which is proportional to the length or travel time of the path and update the pheromone in view of evaporation | Maximal number of iterations or most of ants on a common path | Parameter and the updating mechanism of pheromone |
| SA | The process of annealing in metallurgy | May choose a worse solution with a probability associated with the current temperature | Neighborhood searching around current solution at different temperature | Rise up the temperature gradually and conduct the neighborhood search at each temperature | Maximal temperature or the best solution without change in a series of temperature | The rate of temperature adjustment and the scale of neighborhood search |

## 4. Numerical experiments

To assess the efficiency and effectiveness of the SLLS, we have undertaken several evaluations. Initially, we applied the SLLS to 29 standard benchmark optimization problems and compared its performance against nine well-established heuristic algorithms. Subsequently, we implemented the SLLS to solving the CEC 2014 test suite [60] and compared the results obtained by 8 algorithms. Then, we employed the SLLS to tackle seven renowned engineering optimization problems, conducting a comparative analysis with results obtained from numerous other algorithms. Finally, we delved into investigating the impact of various parameters on the SLLS's performance.

The implementation of the SLLS involved coding it in the JAVA programming language, and the resultant program was executed within the Apache NetBeans IDE 12.6. The computational setup utilized for these experiments featured an Intel(R) Core (TM) i7-1065G7 CPU operating at speeds ranging from 1.30GHz to 1.50GHz, complemented by 16.0GB of RAM.

If no otherwise specified, we will employ the following parameter values for the SLLS. We will set the number of snakes, denoted as $n_S$, to 20, while the maximal number of iterations, $T$, will be fixed at 1000. The truncated range of the natural growth function, represented by $\gamma$, will be assigned a value of 6. For serpentine movement, we will utilize 2 half circles, and for caterpillar movement, 4 touch points will be utilized. The size of the list of visible points will be set to 5.

To compute the initial amplitude of serpentine locomotion, denoted as $L_A(0)$, we will use the following formula:

$$L_A(0) = \frac{1}{5}\sqrt[2]{\sum_{k=1,2,\dots,m} \left(\bar{x}_k - \underline{x}_k\right)^2}. \quad (19)$$

Here, $x_{i,k}$ is constrained within the range $[\underline{x}_k, \overline{x}_k]$ for all $i$ in I, where $k$ ranges from 1 to $m$, as previously mentioned. The value of $m$ corresponds to the dimension of the decision variables. The minimal amplitude of serpentine movement, $L_A^{min}$, will be set to 1.0E-30.

In our experimentation, each problem will be tackled 30 times using the SLLS. In other words, the SLLS program will run independently 30 times for each problem. We will record and compare the average objective values and their associated variations obtained from these 30 runs with results from other algorithms.

**4.1. Standard optimization problems**

In this subsection, we will apply the SLLS to solve 29 standard benchmark optimization problems. These problems encompass seven unimodal problems, six high-dimensional multimodal problems, ten multimodal problems with fixed dimensions, and six composition problems. These benchmark functions have been extensively used to evaluate various heuristic algorithms in previous research. Many well-known algorithms, including MMLA[31], EO [43], GA [3], PSO [4], GSA [34], GWO [18], SSA [21], Evolution Strategy with Covariance Matrix Adaptation (CMA-ES) [61], Success History Based Parameter Adaptation Differential Evolution (SHADE) [62], and SHADE with linear population size reduction hybridized with semi-parameter adaptation of CMA-ES (LSHADE) [63], have been applied to these problems. Detailed descriptions of these problems can be found in Appendix A.

The primary results obtained using the SLLS are summarized in Table 2, alongside results from the other nine algorithms, including MMLA, EO, PSO, GWO, GSA, SSA, CMA-ES, SHADE, and LSHADE.

Among the 29 problems, the SLLS outperforms other algorithms in 17 problems based on the final best average objective values (Ave). We employ Friedman mean rank for the comparison and assign ranks to all algorithms, as shown in Table 2. The top three algorithms are the SLLS, MMLA, and EO, respectively. Notably, the SLLS excels in solving unimodal and multimodal problems but falls behind MMLA when dealing with composition problems. This discrepancy is attributed to MMLA's mechanism of shrinking the search space with varying scales. Given this observation, we believe that incorporating a similar mechanism into the SLLS in the future may enhance its performance in tackling certain composition problems. Additionally, alongside average objective values, Table 2 presents standard deviations (Std), which are consistent with the performance in Ave. In other words, the SLLS demonstrates relatively consistent performance. While the SLLS may not outperform all other algorithms across all problems, its promising performance aligns with the "No Free Lunch Theorem"[2] and underscores the theoretical and practical value of the SLLS.

**Table 2**

Optimization results and comparison for functions.

| Function | SLLS | MMLA | EO | PSO | GWO | GSA | SSA | CMA-ES | SHADE | LSHADE |
|---|---|---|---|---|---|---|---|---|---|---|

| | | | | | | | | | | | | |
|---|---|---|---|---|---|---|---|---|---|---|---|---|
| Unimodal | F1 | Ave | 8.29E-27 | 4.95E−23 | 3.32E−40 | 9.59E−06 | 6.59E−28 | 2.53E−16 | 1.58E−07 | 1.42E−18 | 1.42E−09 | 0.2237 |
| | | Std | 7.12E-27 | 4.92E−24 | 6.78E−40 | 3.35E−05 | 1.58E−28 | 9.67E−17 | 1.71E−07 | 3.13E−18 | 3.09E−09 | 0.1480 |
| | F2 | Ave | 8.43E-7 | 3.12E-10 | 7.12E−23 | 0.02560 | 7.18E−17 | 0.05565 | 2.66293 | 2.98E−07 | 0.0087 | 21.1133 |
| | | Std | 3.35E-7 | 1.97E-11 | 6.36E−23 | 0.04595 | 7.28E−17 | 0.19404 | 1.66802 | 1.7889 | 0.0213 | 9.5781 |
| | F3 | Ave | 1.17E-14 | 2.14E-10 | 8.06E−09 | 82.2687 | 3.29E−06 | 896.534 | 1709.94 | 1.59E−05 | 15.4352 | 88.7746 |
| | | Std | 4.11E-14 | 2.79E-10 | 1.60E−08 | 97.2105 | 1.61E−05 | 318.955 | 11242.3 | 2.21E−05 | 9.9489 | 47.2300 |
| | F4 | Ave | 3.85E-12 | 2.87E-10 | 5.39E−10 | 4.26128 | 5.61E−07 | 7.35487 | 11.6741 | 2.01E−06 | 0.9796 | 2.1170 |
| | | Std | 3.76E-12 | 3.22E-11 | 1.38E−09 | 0.67730 | 1.04E−06 | 1.74145 | 4.1792 | 1.25E−06 | 0.7995 | 0.4928 |
| | F5 | Ave | 8.90E-10 | 24.06265 | 25.32331 | 92.4310 | 26.81258 | 67.5430 | 296.125 | 36.7946 | 24.4743 | 28.8255 |
| | | Std | 3.81E-10 | 0.952351 | 0.169578 | 74.4794 | 0.793246 | 62.2253 | 508.863 | 33.4614 | 11.2080 | 0.8242 |
| | F6 | Ave | 0 | 0 | 8.29E−06 | 8.89E−06 | 0.816579 | 2.5E−16 | 1.80E−07 | 6.83E−19 | 5.31E−10 | 0.2489 |
| | | Std | 0 | 0 | 5.02E−06 | 9.91E−06 | 0.482126 | 1.74E−16 | 3.00E−07 | 6.71E−19 | 6.35E−10 | 0.1131 |
| | F7 | Ave | 2.27E-6 | 0.001634 | 0.001171 | 0.02724 | 0.002213 | 0.08944 | 0.1757 | 0.0275 | 0.0235 | 0.0047 |
| | | Std | 4.13E-6 | 8.85E-04 | 6.54E−04 | 0.00804 | 0.001996 | 0.04339 | 0.0629 | 0.0079 | 0.0088 | 0.0019 |
| Multimodal (High dimensional) | F8 | Ave | -12369.84 | -11863.47 | −9016.34 | −6075.85 | −6123.1 | −2821.1 | −7455.8 | −7007.1 | −11713.1 | −3154.4 |
| | | Std | 103.4623 | 133.4013 | 595.1113 | 754.632 | 909.865 | 493.037 | 772.811 | 773.94 | 230.49 | 317.921 |
| | F9 | Ave | 1.42E-14 | 0 | 0 | 52.8322 | 0.31052 | 25.9684 | 58.3708 | 25.338 | 8.5332 | 67.542 |
| | | Std | 3.55E-15 | 0 | 0 | 16.7068 | 0.35214 | 7.47006 | 20.016 | 8.5539 | 2.1959 | 10.016 |
| | F10 | Ave | 5.02E-14 | 1.28E-12 | 8.34E−14 | 0.00501 | 1.06E−13 | 0.06208 | 2.6796 | 15.587 | 0.3957 | 0.0393 |
| | | Std | 6.68E-14 | 4.90E-13 | 2.53E−14 | 0.01257 | 2.24E−13 | 0.23628 | 0.8275 | 7.9273 | 0.5868 | 0.0151 |
| | F11 | Ave | 2.09E-14 | 4.08E-25 | 0 | 0.02381 | 0.00448 | 27.7015 | 0.0160 | 5.76E−15 | 0.0048 | 0.8948 |
| | | Std | 1.02E-14 | 6.03E-26 | 0 | 0.02870 | 0.00665 | 5.04034 | 0.0112 | 6.18E−15 | 0.0077 | 0.1078 |
| | F12 | Ave | 6.13E-19 | 4.01E-15 | 7.97E−07 | 0.02764 | 0.05343 | 1.79961 | 6.9915 | 2.87E−16 | 0.0346 | 8.18E−04 |
| | | Std | 5.24E-19 | 6.39E-16 | 7.69E−07 | 0.05399 | 0.02073 | 0.95114 | 4.4175 | 5.64E−16 | 0.0875 | 0.0010 |
| | F13 | Ave | 1.92E-17 | 2.44E-24 | 0.029295 | 0.00732 | 0.65446 | 8.89908 | 15.8757 | 3.66E−04 | 7.32E−04 | 0.0102 |
| | | Std | 2.16E-18 | 4.21E-25 | 0.035271 | 0.01050 | 0.00447 | 7.12624 | 16.1462 | 0.0020 | 0.0028 | 0.0103 |
| Multimodal (Fixed-dimensional) | F14 | Ave | 0.998004 | 0.998004 | 0.998004 | 3.84902 | 4.042493 | 5.859838 | 1.1965 | 10.237 | 0.998004 | 1.9416 |
| | | Std | 0 | 0 | 1.54E−16 | 3.24864 | 4.252799 | 3.831299 | 0.5467 | 7.5445 | 5.83E−17 | 2.9633 |
| | F15 | Ave | 0.000307 | 0.000307 | 0.002398 | 0.002434 | 0.00337 | 0.003673 | 0.000886 | 0.0057 | 0.002374 | 3.00E−04 |
| | | Std | 0 | 1.37E-19 | 0.006097 | 0.006081 | 0.00625 | 0.001647 | 0.000257 | 0.0121 | 0.0061 | 1.93E−19 |
| | F16 | Ave | -1.03163 | -1.03163 | −1.03162 | −1.03162 | −1.03163 | −1.03163 | −1.03163 | −1.03162 | −1.03162 | −1.03162 |
| | | Std | 0 | 0 | 6.04E−16 | 6.51E−16 | 2.13E−08 | 4.88E−16 | 6.13E−14 | 6.77E−16 | 6.51E−16 | 1.00E−15 |
| | F17 | Ave | 0.397887 | 0.397887 | 0.397887 | 0.397887 | 0.397889 | 0.397887 | 0.397887 | 0.397887 | 0.397887 | 0.397887 |
| | | Std | 0 | 0 | 0 | 0 | 2.13E−04 | 0 | 3.41E−14 | 0 | 3.24E−16 | 0 |
| | F18 | Ave | 3 | 3 | 3 | 3 | 3.000028 | 3 | 3 | 8.4000 | 3 | 3 |
| | | Std | 0 | 0 | 1.56E−15 | 1.97E−15 | 4.24E−04 | 4.17E−15 | 2.20E−13 | 20.550 | 1.87E−15 | 1.25E−15 |
| | F19 | Ave | -3.86278 | -3.86278 | −3.86278 | −3.86278 | −3.86263 | −3.86278 | −3.86278 | −3.86278 | −3.86278 | −3.86278 |
| | | Std | 0 | 0 | 2.59E−15 | 2.65E−15 | 0.00273 | 2.29E−15 | 1.47E−10 | 2.7E−15 | 2.69E−15 | 2.7E−15 |
| | F20 | Ave | -3.3220 | -3.3220 | −3.2687 | −3.26651 | −3.28654 | −3.31778 | −3.2304 | −3.2903 | −3.27047 | −3.28234 |
| | | Std | 0 | 0 | 0.05701 | 0.06032 | 0.10556 | 0.023081 | 0.0616 | 0.0535 | 0.0599 | 0.0570 |
| | F21 | Ave | -10.1532 | -10.1532 | −8.55481 | −5.9092 | −8.7214 | −5.95512 | −9.6334 | −5.6642 | −9.2343 | −9.4735 |
| | | Std | 0 | 0 | 2.76377 | 3.59559 | 2.6914 | 3.73707 | 1.8104 | 3.3543 | 2.4153 | 1.7626 |
| | F22 | Ave | -10.4029 | -10.4029 | −9.3353 | −7.3360 | −9.2415 | −10.4015 | −9.0295 | −8.4434 | −10.1479 | −10.2258 |
| | | Std | 0 | 0 | 2.43834 | 3.47381 | 1.61254 | 2.01408 | 2.3911 | 3.3388 | 1.3969 | 0.9704 |
| | F23 | Ave | -10.5364 | -10.5364 | −9.63655 | −8.7482 | −10.5343 | −10.5364 | −9.0333 | −8.0750 | −10.2809 | −10.5364 |
| | | Std | 0 | 0 | 2.38811 | 2.55743 | 0.00125 | 2.6E−15 | 2.9645 | 3.5964 | 1.3995 | 1.77E−15 |
| Composition | F24 (CF1) | Ave | 1.02 | 2.41E-27 | 66.666 | 151.18 | 90.229 | 20.000 | 43.333 | 209.48 | 63.333 | 3.3333 |
| | | Std | 0.528 | 7.56E-28 | 95.893 | 123.49 | 105.51 | 48.423 | 67.891 | 215.06 | 80.872 | 18.254 |
| | F25 (CF2) | Ave | 0.04 | 26.34 | 89.837 | 204.92 | 163.56 | 186.77 | 31.133 | 189.83 | 40.508 | 0.0000 |
| | | Std | 0.013 | 22.365 | 56.366 | 118.89 | 89.476 | 62.726 | 52.149 | 170.79 | 61.462 | 0.0000 |
| | F26 (CF3) | Ave | 13.30 | 2.07E-28 | 161.73 | 273.73 | 210.61 | 218.55 | 235.11 | 274.20 | 139.48 | 104.29 |
| | | Std | 4.542 | 1.23E-28 | 33.227 | 110.87 | 95.214 | 117.02 | 80.839 | 213.89 | 33.366 | 14.266 |
| | F27 (CF4) | Ave | 100.16 | 6.31 | 356.44 | 487.45 | 418.63 | 492.33 | 232.44 | 372.99 | 316.62 | 278.63 |
| | | Std | 72.938 | 6.67E-08 | 115.66 | 151.15 | 156.16 | 99.549 | 43.643 | 152.12 | 96.752 | 7.0670 |
| | F28 (CF5) | Ave | 18.67 | 0.60 | 52.309 | 214.56 | 143.81 | 232.32 | 27.538 | 224.85 | 39.515 | 2.02E−17 |
| | | Std | 1.625 | 0.236 | 95.565 | 180.03 | 149.12 | 75.405 | 41.598 | 286.23 | 51.233 | 7.69E−17 |
| | F29 (CF6) | Ave | 400.00 | 160.00 | 768.48 | 794.50 | 837.47 | 845.47 | 628.69 | 845.26 | 684.51 | 540.23 |
| | | Std | 100.09 | 219.09 | 192.94 | 175.94 | 136.45 | 80.524 | 184.48 | 139.52 | 201.22 | 122.75 |
| Friedman mean rank | | | 1.913 | 2.000 | 4.500 | 7.569 | 6.483 | 7.103 | 7.241 | 7.052 | 5.466 | 5.672 |
| Rank | | | 1 | 2 | 3 | 10 | 6 | 8 | 9 | 7 | 4 | 5 |

## 4.2. CEC 2014 functions

In this section, we will evaluate the effectiveness of the SLLS using the well-known CEC 2014 test suite [60]. The results obtained with the SLLS will be compared to those generated by seven other algorithms, namely AHA[25], TLBO[28], GSA, Artificial Bee Colony (ABC) [8], Covariant

Matrix Adaptation with Evolution Strategy (CMA-ES) [64], Success-History Based Adaptive DE (SHADE) [65], and Salp Swarm Algorithm (SSA) [21].

The CEC 2014 test suite comprises 30 problems, including 3 unimodal functions (CF1–CF3), 13 multimodal functions (CF4–CF16), 6 hybrid functions (CF17–CF22), and 8 composition functions (CF23–CF30). To conserve space, we will focus solely on the 30-dimensional versions of these problems. All these problems share the same search range, defined by $[-100,100]^D$, where D=30 represents the number of dimensions. Each problem $\mathrm{CF}i$ is associated with a theoretical optimum such that $F_i^* = 100 \times i$ for all $i \in \{1,2,\ldots,30\}$.

With the exception of the SLLS results, all other data is sourced from literature [25]. To ensure a fair comparison, we maintain a maximum number of function evaluations that does not exceed 25,000, as per the study [25]. For the SLLS, we set $\mathrm{T} = 312$ and $\gamma = 20$, keeping the other parameters consistent with the previous subsection. Specific parameters used by other algorithms can be referenced in Zhao(2022) [25]. Each of the eight algorithms is run independently 30 times, and the results presented here are based on the average performance across these 30 runs.

The outcomes related to the CEC 2014 test suite are summarized in Table 3. Within this table, the best-performing results for each problem are underlined. The data in Table 3 highlights that the SLLS outperforms other algorithms in 12 out of 30 problems. SHADE attains the best results for 7 problems, while AHA leads in 6 other problems. It's worth noting that, except for CMA-ES and SSA, other algorithms also excel in specific problems.

To provide a comprehensive assessment of the overall performance of these eight algorithms, we conducted the Friedman test. The final ranking of each algorithm is presented in Table 3. Notably, the SLLS, AHA, and SHADE secure the top three ranks. The mean rank of the SLLS, which is 2.13, signifies its superior overall performance in tackling these 30 problems.

**Table 3**

Results on functions in CEC 2014.

| Fun. | Index | SLLS | AHA | TLBO | GSA | ABC | CMA-ES | SHADE | SSA |
|---|---|---|---|---|---|---|---|---|---|
| CF1 | Ave | 7.32E+6 | 2.27E+7 | 6.81E+6 | 4.20E+8 | 2.22E+8 | 1.89E+8 | 9.91E+6 | 2.58E+7 |
| | Std | 3.05E+6 | 1.22E+7 | 3.43E+6 | 8.23E+7 | 4.81E+7 | 8.33E+6 | 2.82E+6 | 1.33E+7 |
| CF2 | Ave | 1.54E+4 | 9.42E+6 | 9.85E+5 | 2.33E+10 | 3.07E+8 | 3.89E+9 | 1.01E+4 | 1.03E+4 |
| | Std | 1.15E+4 | 6.98E+6 | 1.16E+6 | 3.26E+9 | 7.11E+7 | 2.27E+9 | 4.40E+3 | 8.60E+3 |
| CF3 | Ave | 1.25E+4 | 6.47E+3 | 4.27E+4 | 8.73E+4 | 8.00E+4 | 2.66E+4 | 5.08E+2 | 6.78E+4 |
| | Std | 9.68E+3 | 3.59E+3 | 6.45E+3 | 6.51E+3 | 1.36E+4 | 1.85E+4 | 1.26E+2 | 1.65E+4 |
| CF4 | Ave | 506.3617 | 5.67E+2 | 5.63E+2 | 2.83E+3 | 7.95E+2 | 1.94E+3 | 5.01E+2 | 5.60E+2 |
| | Std | 33.4790 | 4.77E+1 | 3.69E+1 | 3.37E+2 | 4.52E+1 | 3.64E+2 | 2.92E+1 | 4.12E+1 |
| CF5 | Ave | 520.0196 | 520.1335 | 521.0691 | 519.9992 | 520.0542 | 520.9916 | 521.0091 | 519.9999 |
| | Std | 0.1933 | 9.33E-2 | 4.66E-2 | 4.90E-4 | 5.95E-2 | 8.92E-2 | 4.54E-2 | 1.11E-4 |
| CF6 | Ave | 602.1764 | 6.17E+2 | 6.15E+2 | 6.32E+2 | 6.32E+2 | 6.28E+2 | 6.19E+2 | 6.20E+2 |
| | Std | 1.0022 | 2.68E+0 | 2.46E+0 | 2.74E+0 | 1.37E+0 | 4.93E+0 | 4.90E+0 | 3.33E+0 |
| CF7 | Ave | 700.0018 | 701.0583 | 700.7527 | 956.9402 | 703.4847 | 790.1209 | 700.0042 | 700.0133 |
| | Std | 3.90E-3 | 7.40E-2 | 2.61E-1 | 3.35E+1 | 4.94E-1 | 3.47E+1 | 8.88E-3 | 9.59E-3 |
| CF8 | Ave | 8.03E+2 | 8.48E+2 | 8.67E+2 | 9.45E+2 | 1.02E+3 | 9.17E+2 | 9.10E+2 | 9.42E+2 |
| | Std | 1.5700 | 1.75E+1 | 1.39E+1 | 1.01E+1 | 1.15E+1 | 3.69E+1 | 1.01E+1 | 3.04E+1 |
| CF9 | Ave | 989.6767 | 1.05E+3 | 1.01E+3 | 1.07E+3 | 1.15E+3 | 1.10E+3 | 1.06E+3 | 1.05E+3 |
| | Std | 24.1184 | 3.09E+1 | 3.14E+1 | 1.44E+1 | 1.09E+1 | 4.51E+1 | 8.42E+0 | 4.04E+1 |
| CF10 | Ave | 1.01E+3 | 1.67E+3 | 5.46E+3 | 4.91E+3 | 7.48E+3 | 3.56E+3 | 5.04E+3 | 4.47E+3 |
| | Std | 3.6335 | 3.96E+2 | 1.32E+3 | 3.60E+2 | 2.55E+2 | 9.45E+2 | 5.22E+2 | 6.58E+2 |
| CF11 | Ave | 3.29E+3 | 3.89E+3 | 8.42E+3 | 5.61E+3 | 8.63E+3 | 7.94E+3 | 7.43E+3 | 4.88E+3 |
| | Std | 4.74E+2 | 5.14E+2 | 2.94E+2 | 5.07E+2 | 2.81E+2 | 8.24E+2 | 3.71E+2 | 6.86E+2 |
| CF12 | Ave | 1200.1581 | 1200.2553 | 1203.1427 | 1200.0275 | 1203.0937 | 1202.7605 | 1202.4971 | 1200.6903 |
| | Std | 2.03E-1 | 1.45E-1 | 3.75E-1 | 1.31E-2 | 4.56E-1 | 5.85E-1 | 3.14E-1 | 3.74E-1 |
| CF13 | Ave | 1300.4208 | 1300.4547 | 1300.4339 | 1304.4934 | 1300.6252 | 1302.8761 | 1300.5450 | 1300.4914 |
| | Std | 1.21E-1 | 1.20E-1 | 8.87E-2 | 3.90E-1 | 8.36E-2 | 4.65E-1 | 4.32E-2 | 1.08E-1 |
| CF14 | Ave | 1400.3433 | 1400.3336 | 1400.2649 | 1504.4364 | 1400.2499 | 1445.3904 | 1400.2850 | 1400.3555 |
| | Std | 1.34E-1 | 1.51E-1 | 4.84E-2 | 1.39E+1 | 6.97E-2 | 1.48E+1 | 3.86E-2 | 1.53E-1 |
| CF15 | Ave | 1508.2058 | 1517.8706 | 1527.2681 | 6393.9443 | 1600.5196 | 1602.8043 | 1514.3350 | 1513.4816 |
| | Std | 2.47E+0 | 9.43E+0 | 6.09E+0 | 2.61E+3 | 3.99E+1 | 2.98E+2 | 9.37E-1 | 3.89E+0 |
| CF16 | Ave | 1610.5269 | 1610.9303 | 1612.7586 | 1613.7006 | 1613.0812 | 1612.9198 | 1612.8083 | 1612.4193 |
| | Std | 8.09E-1 | 5.95E-1 | 3.01E-1 | 3.01E-1 | 2.12E-1 | 2.93E-1 | 2.66E-1 | 4.07E-1 |
| CF17 | Ave | 1.45E+6 | 1.32E+6 | 4.29E+5 | 3.77E+7 | 5.36E+6 | 6.70E+6 | 2.88E+5 | 1.57E+6 |
| | Std | 1.32E+6 | 9.18E+5 | 2.67E+5 | 8.06E+6 | 1.42E+6 | 4.99E+6 | 1.72E+5 | 1.41E+6 |
| CF18 | Ave | 2.59E+3 | 3.32E+3 | 2.79E+3 | 1.98E+5 | 1.30E+4 | 1.51E+6 | 1.00E+4 | 6.67E+6 |
| | Std | 2.01E+3 | 1.70E+3 | 1.18E+3 | 1.07E+6 | 1.54E+4 | 5.49E+6 | 6.24E+3 | 4.70E+3 |
| CF19 | Ave | 1911.1944 | 1.93E+3 | 1.92E+3 | 2.10E+3 | 1.92E+3 | 1.97E+3 | 1.91E+3 | 1.92E+3 |
| | Std | 2.08E+1 | 2.83E+1 | 2.41E+1 | 2.36E+1 | 1.31E+0 | 3.03E+1 | 7.83E-1 | 1.59E+1 |
| CF20 | Ave | 2.82E+4 | 1.85E+4 | 2.00E+4 | 2.80E+4 | 4.25E+4 | 3.61E+4 | 3.47E+3 | 2.86E+4 |
| | Std | 1.67E+4 | 7.19E+3 | 5.94E+3 | 1.28E+5 | 1.51E+4 | 2.41E+4 | 2.02E+3 | 1.39E+4 |
| CF21 | Ave | 2.81E+5 | 2.91E+5 | 1.71E+5 | 1.63E+7 | 1.06E+6 | 1.06E+6 | 4.95E+4 | 3.05E+5 |
| | Std | 2.24E+5 | 2.17E+5 | 1.09E+5 | 4.06E+6 | 5.06E+5 | 1.36E+6 | 4.13E+4 | 2.27E+5 |
| CF22 | Ave | 2.64E+3 | 2.85E+3 | 2.60E+3 | 3.70E+3 | 2.89E+3 | 3.04E+3 | 2.65E+3 | 2.67E+3 |
| | Std | 1.21E+2 | 1.94E+2 | 1.40E+2 | 5.65E+2 | 1.41E+2 | 2.64E+2 | 8.48E+1 | 2.20E+2 |
| CF23 | Ave | 2615.2441 | 2500.0000 | 2615.3530 | 2573.1800 | 2621.9003 | 2663.2842 | 2615.2475 | 2633.2248 |
| | Std | 1.0117 | 0.00E+0 | 1.48E-1 | 1.15E+2 | 1.33E+0 | 1.65E+1 | 2.65E-3 | 8.08E+0 |
| CF24 | Ave | 2614.0658 | 2600.0000 | 2600.0800 | 2620.0005 | 2648.2497 | 2636.3027 | 2624.9547 | 2639.0200 |
| | Std | 4.1906 | 0.00E+0 | 1.32E-2 | 5.73E+0 | 3.37E+0 | 6.10E+0 | 9.92E-1 | 8.32E+0 |
| CF25 | Ave | 2704.8905 | 2700.0000 | 2700.0159 | 2705.6276 | 2732.6953 | 2723.8740 | 2711.2547 | 2715.0444 |
| | Std | 1.1362 | 0.00E+0 | 8.70E-2 | 2.69E+0 | 3.65E+0 | 6.90E+0 | 1.29E+0 | 3.05E+0 |
| CF26 | Ave | 2700.2887 | 2763.5233 | 2713.7514 | 2793.1042 | 2707.7046 | 2703.2826 | 2700.8536 | 2700.5293 |
| | Std | 1.05E-1 | 4.88E+1 | 3.44E+1 | 2.16E+1 | 1.96E+1 | 1.08E+0 | 6.37E-2 | 1.46E-1 |
| CF27 | Ave | 3058.7903 | 2900.0000 | 3271.3167 | 4557.5228 | 3502.8503 | 3470.9875 | 3108.8516 | 3194.4532 |
| | Std | 12.2958 | 0.00E+0 | 1.25E+2 | 3.96E+2 | 1.22E+2 | 2.29E+2 | 8.07E+1 | 1.64E+2 |
| CF28 | Ave | 3157.1214 | 3000.0000 | 3938.6520 | 5215.9261 | 4357.7299 | 7445.9885 | 3929.0812 | 4653.9440 |
| | Std | 61.7996 | 0.00E+0 | 1.42E+2 | 8.98E+2 | 1.18E+2 | 5.62E+2 | 7.64E+1 | 5.52E+2 |
| CF29 | Ave | 3.22E+3 | 3.10E+3 | 4.69E+3 | 6.25E+6 | 3.78E+5 | 2.85E+7 | 1.26E+4 | 1.35E+4 |
| | Std | 2.51E+3 | 0.00E+0 | 9.96E+2 | 2.02E+7 | 1.97E+5 | 3.16E+7 | 2.41E+3 | 5.47E+3 |
| CF30 | Ave | 3.87E+3 | 4.45E+3 | 6.98E+3 | 2.57E+6 | 4.70E+4 | 1.01E+5 | 8.44E+3 | 6.39E+4 |
| | Std | 2.66E+2 | 3.53E+3 | 1.52E+3 | 6.84E+5 | 1.32E+4 | 6.05E+4 | 9.29E+2 | 6.13E+4 |
| Friedman mean rank | | 2.13 | 3.05 | 3.80 | 6.45 | 6.27 | 6.42 | 3.37 | 4.48 |
| Rank | | 1 | 2 | 4 | 8 | 6 | 7 | 3 | 5 |

## 4.3. Engineering optimization problems

In this subsection, we will apply the SLLS to solve seven well-known engineering optimization problems. These problems are derived from real-life scenarios and pose a considerable challenge to heuristic algorithms due to their complex search spaces and multiple local minima. Table 4 provides an overview of the key characteristics of these engineering optimization problems. The "Mixed/Con/Discrete" column in Table 4 indicates the nature of the variables involved. "Discrete" signifies that all variables are discrete, "continuous" implies that all variables fall within corresponding continuous intervals, and "Mixed" indicates a combination of discrete and continuous variables. Detailed mathematical formulations of these seven engineering problems are provided in Appendix B. It's worth noting that we modified the original model of the Robot gripper optimization problem by adding three constraints to ensure result feasibility.

To facilitate meaningful comparisons with other algorithms in terms of the total Number of Objective Function Evaluations (NFE), we will adjust the total number of iterations of the SLLS to align the resulting NFEs with those of other algorithms.

In this paper, we employ penalty terms to address constraint violations. When a constraint is breached, the absolute difference between the boundary and the actual function value resulting from the constraint is multiplied by a large positive constant and added to the objective function.

Table 5 presents a comprehensive summary of the primary results obtained for these engineering problems using the SLLS.

**Table 4**

The main characters of the seven engineering optimization problems.

| No. | Name of problem | Number of Variables | Number of Constraints | Mixed/Con/Discrete |
|---|---|---|---|---|
| 1 | Multiple disc clutch brake design | 5 | 8 | Discrete |
| 2 | Robot gripper | 7 | 10 | Continuous |
| 3 | Rolling element bearing | 10 | 7 | Mixed |
| 4 | Hydrodynamic thrust bearing | 4 | 7 | Continuous |
| 5 | Belleville spring | 4 | 7 | Continuous |
| 6 | Step-cone pulley | 5 | 8 | Continuous |
| 7 | Speed reducer design | 7 | 11 | Continuous |

**Table 5**

The main results from SLLS for the seven engineering problems.

| No. of Problem | 1 | 2 | 3 | 4 | 5 | 6 | 7 |
|---|---|---|---|---|---|---|---|
| Best | 0.313657 | 4.91124796 | 81859.7415 | 1625.443 | 1.9807 | 19.1331 | 2994.5735 |
| Mean | 0.313657 | 5.00612329 | 81858.4177 | 1627.247 | 1.9930 | 21.3302 | 2996.7624 |
| Worst | 0.313657 | 6.06680250 | 81793.2266 | 1631.863 | 2.0118 | 22.6464 | 2998.2019 |
| Std | 0.0 | 0.241 | 5.269 | 2.381 | 0.027 | 1.013 | 1.822 |
| NFE | 1200 | 36000 | 16000 | 48000 | 24000 | 72000 | 32000 |

To demonstrate the effectiveness of the SLLS, we have compiled Tables 6 to 12, which compare the results obtained using the SLLS with those from other algorithms. In order to ensure a fair comparison, the results from other algorithms have been sourced from relevant literature and represent their best-performing outcomes.

Following the sequence number provided in Table 4, the first problem has been addressed by various algorithms, including AHA [25], TLBO [28], ABC [28], EOBL-GOA [46], AVOA [24], AEO [44], HHO [23], FSO [22], PVS [41], and WCA [36]. The second problem has been tackled using multiple existing heuristic methods, such as an approach based on a combination of the grasshopper optimization algorithm and the Nelder–Mead algorithm (HGOANM) [66], TLBO [28], PVS [41] and a way using Force and Displacement Transmission Ratio(2014)[67]. The third problem has been previously explored using EAPSO [47], TLBO [28], AHA [25], AVOA [24], AEO [44], HHO [23], MBA [39] and WCA [36]. The fourth problem has been optimized using AHA [25], TLBO [28], EOBL-GOA [46], and PVS [41]. The fifth problem has been addressed with the involvement of AHA [25], PVS [41], MBA [39], TLBO [28], a multi-objective optimization approach (Coello) [68], the combined genetic search technique (Gene AS I and Gene AS II) [69], and the APPROX (Griffith and Stewart's successive linear approximation) [70]. The sixth problem has been studied by researchers using PVS [41], ABC [28], TLBO [28], and an improved constrained Differential Evolution variant (rank-iMDDE) [71]. The seventh problem has been optimized with AHA [25], Geometric Mean Optimizer(GMO)[49], EAPSO [47], AEO [44], PVS [41], MBA [39], CS [72], BA [17], WCA [36], and Artificial Bee Colony (ABC) algorithm[73].

The best objective values for each problem are highlighted in red. Similarly, the best means of the objective values, representing the average values of 30 independent runs, and the associated standard deviations are also highlighted in red. The comparison reveals that the SLLS achieves the best objective values in six out of the seven problems, showcasing its promise as a heuristic algorithm. Only in the context of the speed reducer design optimization does the SLLS fall short of the MMLA's performance. However, when compared to other algorithms, the SLLS maintains a highly promising level of performance.

It is essential to clarify that the feasibility of the results presented in Tables 6 to 12 has been thoroughly verified. Some results may exhibit slight violations of the related constraints. When such violations occur, the corresponding data entries are underlined in the tables. Two factors may account for these violations. Firstly, some related studies may have employed low-accuracy criteria. Secondly, the results presented in certain papers may have been truncated to meet the required formatting constraints.

**Table 6**

Results for the multiple disc clutch brake design.

| | | SLLS | MMLA | AHA | TLBO | EOBL-GOA | AEO | HHO | FSO |
|---|---|---|---|---|---|---|---|---|---|
| $x$ | $x_1$ | 70 | 70 | 70 | 70 | 70 | 70 | 70 | 70 |
| | $x_2$ | 90 | 90 | 90 | 90 | 90 | 90 | 90 | 90 |
| | $x_3$ | 1 | 1 | 1 | 1 | 1 | 1 | 1 | 1 |
| | $x_4$ | 810 | 830 | 840 | 810 | 984 | 810 | 1000 | 870 |
| | $x_5$ | 3 | 3 | 3 | 3 | 3 | 3 | 2.313 | 3 |
| $g(x)$ | $g_1(x)$ | 0 | 0 | 0 | 0 | 0 | 0 | 0 | 0 |
| | $g_2(x)$ | 24 | 24 | 24 | 24 | 24 | 24 | 25.5 | 24 |
| | $g_3(x)$ | 0.919 | 0.917 | 0.916 | 0.919 | 0.902 | 0.919 | 0.901 | 0.913 |
| | $g_4(x)$ | 9.830 | 9.826 | 9.824 | 9.830 | 9.794 | 9.830 | 9.791 | 9.818 |
| | $g_5(x)$ | 7.895 | 7.895 | 7.895 | 7.895 | 7.895 | 7.895 | 7.895 | 7.895 |
| | $g_6(x)$ | 0.263 | 0.619 | 1.198 | 0.702 | 2.869 | 0.702 | -2.905 | 0.876 |
| | $g_7(x)$ | 37.706 | 40.119 | 41.325 | 37.706 | 58.695 | 37.706 | 20.417 | 44.944 |
| | $g_8(x)$ | 14.737 | 14.381 | 13.802 | 14.298 | 12.131 | 14.298 | 17.905 | 14.124 |
| Best | | 0.313657 | 0.313657 | 0.313657 | 0.313657 | 0.313657 | 0.313657 | 0.259769 | 0.313657 |
| Mean | | 0.313657 | 0.315225 | 0.321684 | 0.327166 | 0.313657 | 0.321684 | NA | 0.313942 |
| Worst | | 0.313657 | 0.321498 | 0.333260 | 0.392071 | 0.321612 | 0.333260 | NA | 0.333705 |
| std | | 0.0 | 0.003137 | 0.008 | 0.67 | 0.264 | 0.008154 | NA | 4.25e-6 |
| NFE | | 1200 | 600 | 600 | 600 | 300 | 500 | NA | 400 |

**Table 7**

Results for the robot gripper design.

| | | SLLS | HGOANM | TLBO | PVS |
|---|---|---|---|---|---|
| $x$ | $x_1$ | 149.954 | 150.000 | 150 | 150 |
| | $x_2$ | 119.441 | 149.883 | 150 | 150 |
| | $x_3$ | 200.0 | 200 | 200 | 200 |
| | $x_4$ | 27.535 | 0 | 0 | 0 |
| | $x_5$ | 145.927 | 150 | 150 | 150 |
| | $x_6$ | 158.477 | 100.943 | 100 | 100 |
| | $x_7$ | 2.743 | 2.297 | 2.340 | 2.312 |
| $g(x)$ | $g_1(x)$ | 39.003 | 49.989 | -537.514 | -545.098 |
| | $g_2(x)$ | 10.997 | 0.011 | 587.514 | 595.098 |
| | $g_3(x)$ | 47.218 | 49.992 | 33.650 | 43.777 |
| | $g_4(x)$ | 2.782 | 0.008 | 16.350 | 6.223 |
| | $g_5(x)$ | 46700.191 | 79740.197 | 80000.0 | 80000.0 |
| | $g_6(x)$ | 4139.838 | 36.009 | 0.0 | 0.0 |
| | $g_7(x)$ | 58.477 | 0.943 | 0.0 | 0.0 |
| | $g_8(x)$ | 34.123 | 0.826 | 0.0 | 0.0 |
| | $g_9(x)$ | 130.339 | 100.826 | 100.0 | 100.0 |
| | $g_{10}(x)$ | 108.543 | 198.940 | 200.0 | 200.0 |
| Best | | 4.91124796 | 71.05999 | 74.9999999 | 74.99999 |
| Mean | | 5.00612329 | NA | NA | NA |
| Worst | | 6.06680250 | NA | NA | NA |
| std | | 0.241 | NA | NA | NA |
| NFE | | 36000 | NA | 25000 | 25000 |

**Table 8**

Results for the rolling element bearing design.

| | | SLLS | MMLA | AHA | EAPSO | AVOA | AEO | HHO | MBA |
|---|---|---|---|---|---|---|---|---|---|
| $x$ | $x_1$ | 125.719 | 125.719 | 125.718 | 125.719 | 125.723 | 125.719 | 125.0 | 125.715 |
| | $x_2$ | 21.426 | 21.426 | 21.425 | 21.426 | 21.423 | 21.426 | 21.0 | 21.423 |
| | $x_3$ | 10.749 | 10.920 | 10.528 | 11 | 11.001 | 11.395 | 11.092 | 11.0 |
| | $x_4$ | 0.515 | 0.515 | 0.515 | 0.515 | 0.515 | 0.515 | 0515 | 0.515 |
| | $x_5$ | 0.515 | 0.515 | 0.515 | 0.515 | 0.515 | 0.515 | 0.515 | 0.515 |
| | $x_6$ | 0.404 | 0.479 | 0.470 | 0.5 | 0.404 | 0.410 | 0.4 | 0.489 |
| | $x_7$ | 0.700 | 0.667 | 0.641 | 0.644 | 0.619 | 0.638 | 0.6 | 0.628 |
| | $x_8$ | 0.300 | 0.300 | 0.300 | 0.3 | 0.3 | 0.300 | 0.3 | 0.300 |
| | $x_9$ | 0.020 | 0.081 | 0.095 | 0.047 | 0.069 | 0.047 | 0.050 | 0.097 |
| | $x_{10}$ | 0.600 | 0.629 | 0.682 | 0.605 | 0.602 | 0.670 | 0.6 | 0.646 |
| $g(x)$ | $g_1(x)$ | 3.304E-13 | 3.066e-10 | 3.900E-5 | -8.657E-10 | 0.001 | 3.148E-8 | 0.106 | 5.638E-4 |
| | $g_2(x)$ | 14.576 | 9.328 | 9.936 | 7.851 | 14.537 | 14.137 | 14.0 | 8.630 |
| | $g_3(x)$ | 6.149 | 3.806 | 2.007 | 2.199 | 0.461 | 1.838 | 0.0 | 1.101 |
| | $g_4(x)$ | -3.426 | -2.551 | -0.958 | -3.265 | -3.349 | -3.141 | -3.0 | -2.041 |
| | $g_5(x)$ | 0.719 | 0.719 | 0.718 | 0.719 | 0.723 | 0.719 | 0.0 | 0.715 |
| | $g_6(x)$ | 4.281 | 19.499 | 23.062 | 11.086 | 16.560 | 11.023 | 12.619 | 23.611 |
| | $g_7(x)$ | 6.307 | 3.281E-7 | 2.574E-4 | 8.000E-9 | 6.300E-6 | 2.548 | 0.7 | 5.179E-4 |
| | $g_8(x)$ | 0.0 | 3.667E-12 | 0.0 | 0.0 | 0.0 | 0.0 | 0.0 | 0 |
| | $g_9(x)$ | 0.0 | 1.775E-11 | 1.550E-4 | 0.0 | 0.0 | 9.000E-7 | 0.0 | 0 |
| Best | | 81859.7415 | 81859.7392 | 81812.0128 | 81859.7416 | 81843.4096 | 81859.2991 | 78897.8107 | 81843.686 |
| Mean | | 81858.4177 | 81777.3303 | NA | NA | NA | NA | NA | NA |
| Worst | | 81793.2266 | 81344.6726 | NA | NA | NA | NA | NA | NA |
| std | | 5.269 | 120.4525 | NA | NA | NA | NA | NA | NA |
| NFE | | 16000 | 15000 | 15000 | 15000 | NA | 10000 | NA | 15100 |

**Table 9**

Results for the hydrodynamic thrush bearing design.

| | | SLLS | MMLA | TLBO | AHA | EOBL-GOA | PVS |
|---|---|---|---|---|---|---|---|
| $x$ | $x_1$ | 5.955 | 5.955 | 5.956 | 5.956 | 5.956 | 5.956 |
| | $x_2$ | 5.389 | 5.389 | 5.389 | 5.389 | 5.389 | 5.389 |
| | $x_3$ | 5.358E-6 | 5.359E-6 | 5.4E-6 | 5.359E-6 | 5.359E-6 | 1.0E-5 |
| | $x_4$ | 2.269 | 2.269 | 2.270 | 2.270 | 2.26966338 | 2.270 |
| $g(x)$ | $g_1(x)$ | 2.917E-4 | 1.713E-8 | 1.375E-4 | 0.041 | NA | NA |
| | $g_2(x)$ | 8.112E-5 | 4.038E-9 | 1.0E-6 | 1.811E-4 | NA | NA |
| | $g_3(x)$ | 3.583E-4 | 8.96132e-7 | 2.0E-8 | 7.629E-4 | NA | NA |
| | $g_4(x)$ | 3.244E-4 | 3.24363e-4 | 3.244E-4 | 3.244E-4 | NA | NA |
| | $g_5(x)$ | 0.567 | 0.567 | 0.567 | 0.567 | NA | NA |
| | $g_6(x)$ | 8.334E-4 | 8.334E-4 | 8.334E-4 | 9.964E-4 | NA | NA |
| | $g_7(x)$ | 0.004 | 2.523E-9 | 9.08E-6 | 4.646e-3 | NA | NA |
| Best | | 1625.443 | 1625.443 | 1625.443 | 1625.450 | 1625.443 | 1625.444 |
| Mean | | 1627.247 | 1644.001 | 1797.708 | 1680.781 | 1645.422 | 1832.492 |
| Worst | | 1631.863 | 1935.610 | 2096.801 | 1850.381 | 1675.493 | NA |
| std | | 2.381 | 61.108 | NA | 57 | 12.615 | NA |
| NFE | | 48000 | 50000 | 50000 | 50000 | 20000 | 25000 |

**Table 10**

Results for the Belleville spring design.

| | | SLLS | MMLA | AHA | TLBO | PVS | MBA | Coello | Gene ASII |
|---|---|---|---|---|---|---|---|---|---|
| $x$ | $x_1$ | 0.204 | 0.204 | 0.204 | 0.204 | 0.204 | 0.204 | 0.208 | 0.210 |
| | $x_2$ | 0.200 | 0.200 | 0.200 | 0.2 | 0.2 | 0.2 | 0.2 | 0.204 |
| | $x_3$ | 10.025 | 10.012 | 10.030 | 10.030 | 10.030 | 10.030 | 8.751 | 9.268 |
| | $x_4$ | 12.006 | 11.995 | 12.010 | 12.01 | 12.01 | 12.01 | 11.067 | 11.499 |
| $g(x)$ | $g_1(x)$ | 1.933 | 1.262E-5 | 2.235 | 0.506 | 2.498 | 0.243 | 2145.411 | 2127.262 |
| | $g_2(x)$ | 2.057 | 6.285E-6 | -0.267 | -0.036 | -0.274 | -0.028 | 39.7502 | 194.223 |
| | $g_3(x)$ | 1.607E-5 | 4.715E-10 | 0.0 | 0.0 | 0.0 | 0.0 | 0.0 | 0.004 |
| | $g_4(x)$ | 1.596 | 1.596 | 1.596 | 1.596 | 1.596 | 1.596 | 1.592 | 1.586 |
| | $g_5(x)$ | 0.004 | 0.015 | 5.000E-6 | 0.0 | 0.0 | 0.0 | 0.943 | 0.511 |
| | $g_6(x)$ | 1.981 | 1.983 | 1.980 | 1.980 | 1.980 | 1.980 | 2.316 | 2.231 |
| | $g_7(x)$ | 0.199 | 0.199 | 0.199 | 0.199 | 0.199 | 0.199 | 0.214 | 0.209 |
| Best | | 1.9807 | 1.9807 | 1.9797 | 1.9797 | 1.9797 | 1.9797 | 2.1220 | 2.1626 |
| Mean | | 1.9930 | 2.0437 | 1.9860209 | 1.9886513 | 1.983524 | 1.984698 | NA | NA |
| worst | | 2.0118 | 2.1910 | 2.1041923 | 2.1149819 | NA | 2.005431 | NA | NA |
| std | | 0.027 | 0.054 | 0.023 | 0.027 | NA | 7.78E−03 | NA | NA |
| NFE | | 24000 | 25000 | 24000 | 24000 | 15000 | 10600 | 24000 | 24000 |

**Table 11**

Results for the step-cone pulley design.

| | | SLLS | MMLA（Strict） | TLBO | PVS | ABC | Rank-iMDDE |
|---|---|---|---|---|---|---|---|
| $x$ | $x_1$ | 35.871 | 36.066 | 40 | 40 | NA | 100 |
| | $x_2$ | 49.357 | 49.625 | 54.764 | 54.764 | NA | 34.582 |
| | $x_3$ | 65.804 | 66.162 | 73.013 | 73.013 | NA | 47.582 |
| | $x_4$ | 78.903 | 79.332 | 88.428 | 88.428 | NA | 63.438 |
| | $x_5$ | 96.397 | 95.983 | 85.986 | 85.986 | NA | 76.067 |
| $g(x)$ | $h_1(x)$ | 0.0 | 0.0 | 0.001 | 0.001 | NA | 0.371 |
| | $h_2(x)$ | 0.0 | 0.0 | 9.999E-4 | 9.999E-4 | NA | 0.367 |
| | $h_3(x)$ | 0.0 | 0.0 | 0.001 | 0.001 | NA | 0.352 |
| | $g_1(x)$ | 0.989 | 0.988 | 0.987 | 0.986864 | NA | 0.963 |
| | $g_2(x)$ | 0.999 | 0.998 | 0.997 | 0.997360 | NA | 0.999 |
| | $g_3(x)$ | 1.009 | 1.009 | 1.011 | 1.010154 | NA | 1.008 |
| | $g_4(x)$ | 1.019 | 1.019 | 1.021 | 1.020592 | NA | 1.0156 |
| | $g_5(x)$ | 705.658 | 707.061 | 698.577 | 698.577228 | NA | 2211.360 |
| | $g_6(x)$ | 486.668 | 487.841 | 475.827 | 475.827090 | NA | 19.148 |
| | $g_7(x)$ | 216.927 | 217.805 | 209.037 | 209.036940 | NA | -116.522 |
| | $g_8(x)$ | 5.804E-6 | 0.635 | -1.528 | -1.528 | NA | -204.632 |
| Best | | 19.1331 | 19.2588 | 21.1797 | 21.1797 | 16.6347 | 36.8448 |
| Mean | | 21.3302 | 24.7556 | 24.0114 | NA | 36.0995 | NA |
| worst | | 22.6464 | 29.1133 | 74.0230 | NA | 145.4705 | NA |
| std | | 1.013 | 2.187 | 0.34 | NA | 0.06 | NA |
| NFE | | 72000 | 75000 | NA | 15000 | 15000 | NA |

**Table 12**

Results for the speed reducer design.

| | | SLLS | MMLA | GMO | AHA | EAPSO | AEO | MBA | CS | BA | WCA |
|---|---|---|---|---|---|---|---|---|---|---|---|
| $x$ | $x_1$ | 3.5 | 3.5 | 3.5 | 3.5 | 3.5 | 3.5 | 3.5 | 3.502 | 3.5 | 3.5 |
| | $x_2$ | 0.7 | 0.7 | 0.7 | 0.7 | 0.7 | 0.7 | 0.7 | 0.7 | 0.7 | 0.7 |
| | $x_3$ | 17.0 | 17.0 | 17.0 | 17.0 | 17.0 | 17 | 17 | 17 | 17 | 17 |
| | $x_4$ | 7.3 | 7.3 | 7.3 | 7.300 | 7.3 | 7.3 | 7.300 | 7.605 | 7.300 | 7.3 |
| | $x_5$ | 7.716 | 7.715 | 7.728 | 7.715 | 7.715 | 7.7153 | 7.716 | 7.818 | 7.715 | 7.715 |
| | $x_6$ | 3.351 | 3.350 | 3.350 | 3.350 | 3.350 | 3.3502 | 3.350 | 3.352 | 3.350 | 3.350 |
| | $x_7$ | 5.287 | 5.287 | 5.287 | 5.287 | 5.287 | 5.287 | 5.287 | 5.288 | 5.288 | 5.287 |
| $g(x)$ | $g_1(x)$ | -0.074 | -0.074 | -0.074 | -0.074 | -0.074 | -0.074 | -0.074 | -0.074 | -0.074 | -0.739 |
| | $g_2(x)$ | -0.198 | -0.198 | -0.198 | -0.198 | -0.198 | -0.198 | -0.198 | -0.198 | -0.198 | -0.198 |
| | $g_3(x)$ | -0.499 | -0.499 | -0.499 | -0.499 | -0.499 | -0.499 | -0.499 | -0.435 | -0.499 | -0.499 |
| | $g_4(x)$ | -0.905 | -0.905 | -0.904 | -0.905 | -0.905 | -0.905 | -0.905 | -0.901 | -0.905 | -0.905 |
| | $g_5(x)$ | -3.130E-4 | -1.506E-10 | 1.313E-5 | -2.409E-8 | -3.496E-9 | 1.432E-13 | -2.930E-6 | -0.001 | 4.195E-6 | 5.965E-7 |
| | $g_6(x)$ | -4.836E-6 | -7.861E-11 | -2.344E-5 | -1.036E-9 | 2.824E-9 | 5.918E-13 | 3.505E-7 | -4.0E-4 | -4.797E-4 | 2.637E-7 |
| | $g_7(x)$ | -0.703 | -0.702 | -0.703 | -0.703 | -0.703 | -0.703 | -0.703 | -0.703 | -0.703 | -0.703 |
| | $g_8(x)$ | 0.0 | -7.372E-11 | 0.0 | -2.431E-8 | 0.0 | 6.8E-13 | 0.0 | -4.0E-4 | 0.0 | 0.0 |
| | $g_9(x)$ | -0.583 | -0.583 | -0.583 | -0.796 | -0.583 | -0.796 | -0.583 | -0.583 | -0.583 | -0.5833 |
| | $g_{10}(x)$ | -0.051 | -0.051 | -0.051 | -0.051 | -0.051 | -0.051 | -0.051 | -0.089 | -0.051 | -0.0513 |
| | $g_{11}(x)$ | -4.548E-5 | -1.434E-13 | -0.002 | -1.776E-8 | 7.777E-10 | 9.015E-14 | -5.866E-5 | -0.013 | 1.205E-4 | 5.185E-8 |
| Best | | 2994.5735 | 2994.4711 | 2994.750 | 2994.4716 | 2994.471 | 2994.4711 | 2994.4825 | 3000.981 | 2994.4671 | 2994.4711 |
| Mean | | 2996.7624 | 2994.4711 | NA | 2994.4717 | NA | 2994.4711 | 2999.6524 | 3007.1997 | 2994.4671 | 2994.4744 |
| worst | | 2998.2019 | 2994.4717 | NA | 2994.4732 | NA | 2994.4711 | 2999.6524 | 3009.9 | 4973.8644 | 2994.5056 |
| std | | 1.822 | 2.976E-7 | NA | 4.251E-4 | NA | 1.239E-7 | 1.56 | 4.963 | 721.518 | 7.4E-3 |
| NFE | | 32000 | 26000 | NA | 30000 | 15000 | 22000 | 6300 | 250000 | 15000 | 15150 |

## 4.4. Sensitivity analysis of the main parameters

As is well-known, the performance of a heuristic algorithm is often influenced by the selection of its parameters. Therefore, it becomes necessary to investigate the impact of parameter changes on the performance of a heuristic algorithm in order to guide the proper parameter selection.

For the SLLS, there are seven key parameters that merit investigation. In this section, we will conduct a sensitivity analysis of these parameters.

To carry out this analysis, we will utilize the well-known Weierstrass function as the target function. This choice is motivated by the function's complexity, particularly when considering its high-dimensional version. The Weierstrass function is represented as follows:

$$f(x) = \sum_{i=1}^{n}\left(\sum_{k=0}^{k_{max}}\left[a^k \cos\left(2\pi b^k (x_i + 0.5)\right)\right]\right) - n\sum_{k=0}^{k_{max}}\left[a^k \cos(\pi b^k)\right], \tag{20}$$

where $a = 0.5$, $\mathrm{b} = 3$, $k_{max} = 20$, and $x \in [-0.5, 0.5]^n$. $n$ is the dimension of this function. This function has many local minima but only one global minimum at $x^* = (0, \dots, 0)$, $f(x^*) = 0$. Figure 8 illustrates the two-dimensional version of this function. However, it's important to mention that our upcoming experiments will focus on the five-dimensional Weierstrass function. As we've done previously, we'll execute the program 30 times to tackle this problem for each scenario linked to a specific parameter value.

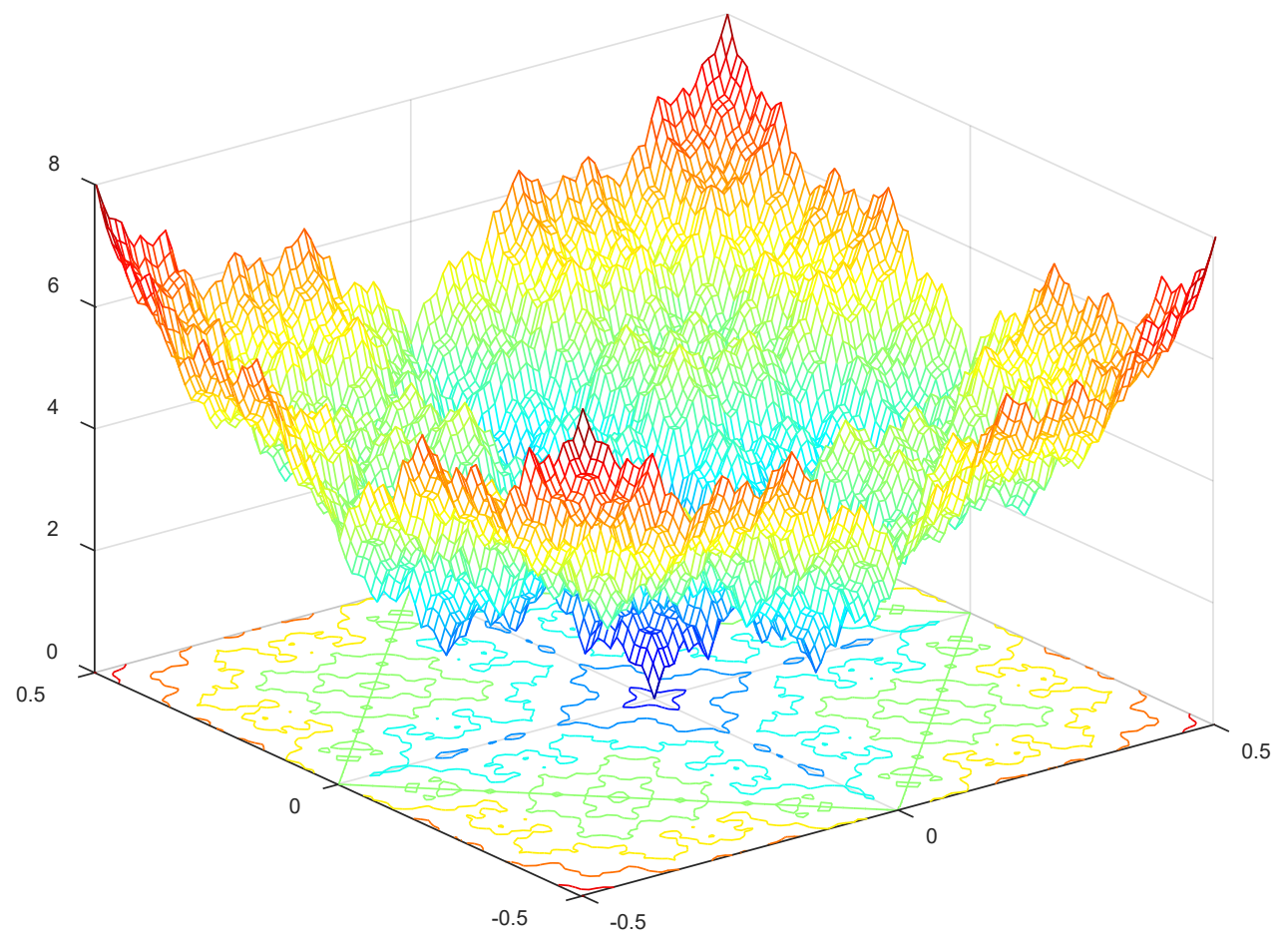

Fig. 8. The 2D version of Weierstrass Function.

The seven parameters to be analyzed are as follows: the number of snakes in the swarm, the radius of the truncated range of the natural growth function $\gamma$, the maximal number of iterations $T$, the size of the list of visible points, the number of half circles in a serpentine locomotion, the number of touch points in the caterpillar locomotion, and the rate for demarcation $r_{CL}$. The benchmark values for these parameters, in the order listed above, are 20, 6, 500, 5, 2, 4, and 0.5, respectively.

Before delving into the sensitivity analysis of parameters, let's first explore what happens to the visible spots and an average snake during a typical search process. With the benchmark values of parameters mentioned above, we ran the program once to solve the five-dimensional Weierstrass function, resulting in a final best fitness value of 5.27E-4. This value corresponds to specific elements: 6.45E-8, -7.99E-11, -5.42E-8, 4.69E-8, and 4.75E-8. The computation took 741 milliseconds.

Throughout the 500 iterations of the search process, the ordered list of visible spots is continuously updated. Figure 9 illustrates the changing objective values of these five visible spots. In the figure, each curve associated with the $i$th visible spot is represented as vs(i). These curves consistently display a common decreasing trend as the number of iterations increases.

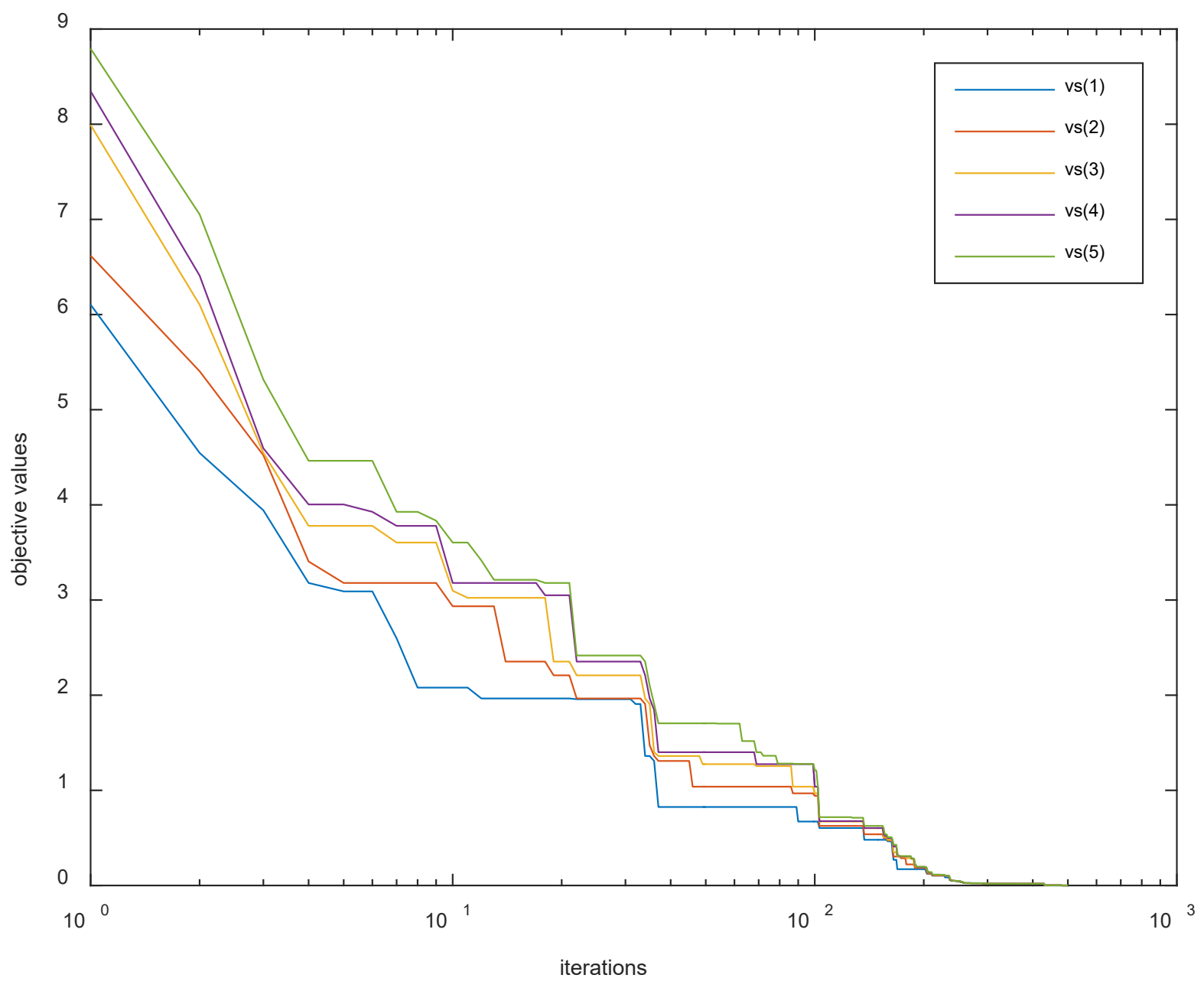

Fig. 9. The variation of objective values of visible spots.

In contrast to the visible spots, which consistently exhibit a continuous decrease in objective values, the trail followed by a snake is marked by a series of spots with fluctuating objective values. Figure 10 portrays the objective values of touch spots in a typical trail taken by a single snake. Given that there are four touch spots within one movement (corresponding to one iteration of the SLLS), a snake's trail encompasses 2000 touch spots.

On one hand, the oscillating objective values of these touch spots along the trail reflect the random locomotion choices made by the snake. On the other hand, these changes signify the intricate variations in the Weierstrass function's values across its entire feasible search space.

Fig. 10. The changing objective values for a snake in the searching process.

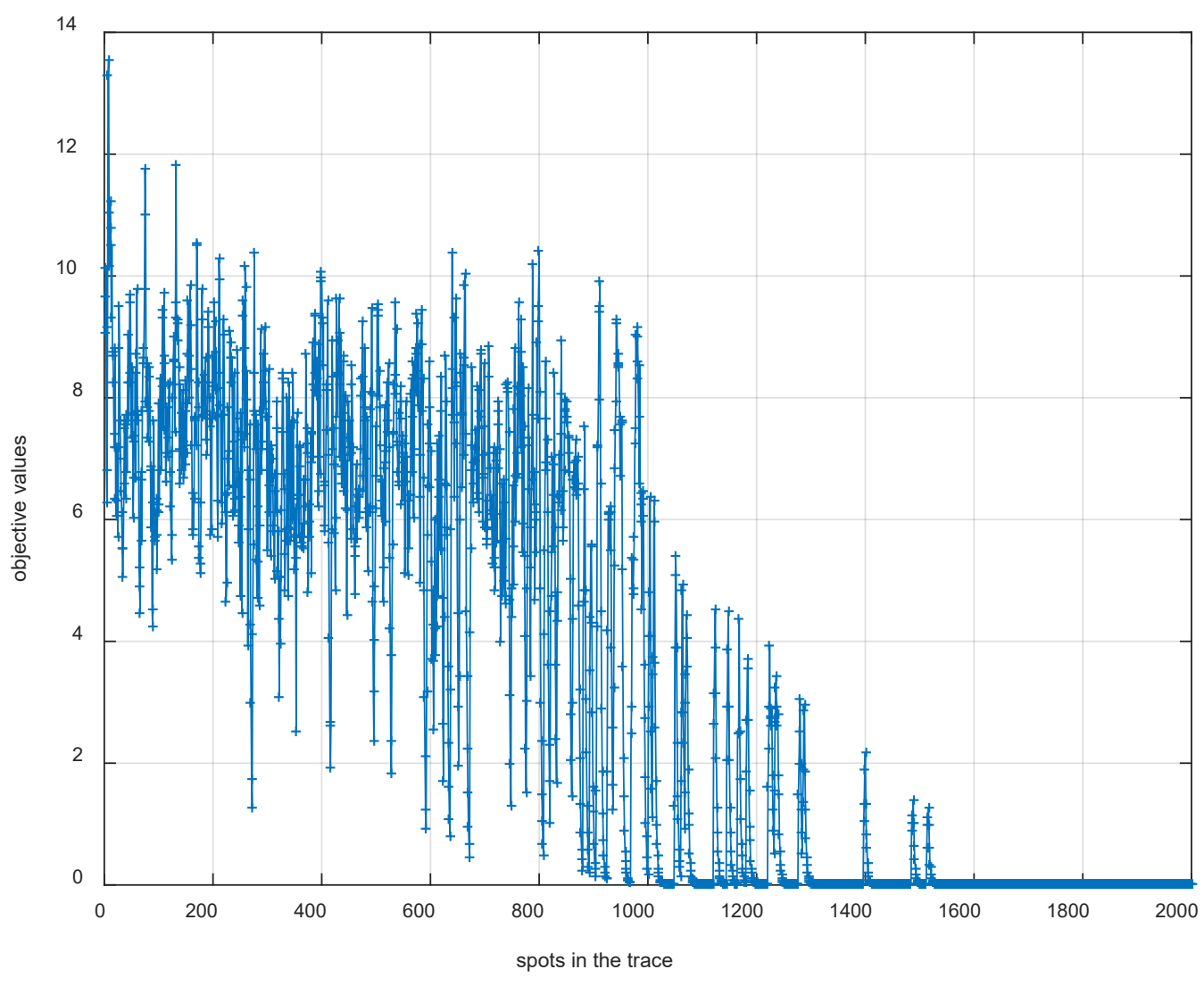

In correspondence with Figure 10, Figure 11 illustrates the selection of locomotion modes during the search process. In Figure 11, a point located at (t, -1) marked with an asterisk (*) signifies that a serpentine locomotion is employed in iteration t, while a point located at (t, 1) indicates the choice of caterpillar locomotion in iteration t. The consecutive choices of locomotion are connected by solid lines.

The evolving preferences in locomotion displayed in Figure 11 for a typical snake reveal that, during the initial phase of the search process, the snake typically favors exploring the search space using serpentine locomotion. However, as the number of iterations increases, the snake transitions towards exploiting previously examined promising locations with caterpillar locomotion. In the concluding phase, the focus of the snake is primarily on exploitation. This behavior mirrors the adaptive influence of the learning efficiency generated by the Sigmoid function. Over the course of 500 iterations, the trace includes a total of 265 instances of caterpillar locomotion and 235 instances of serpentine locomotion.

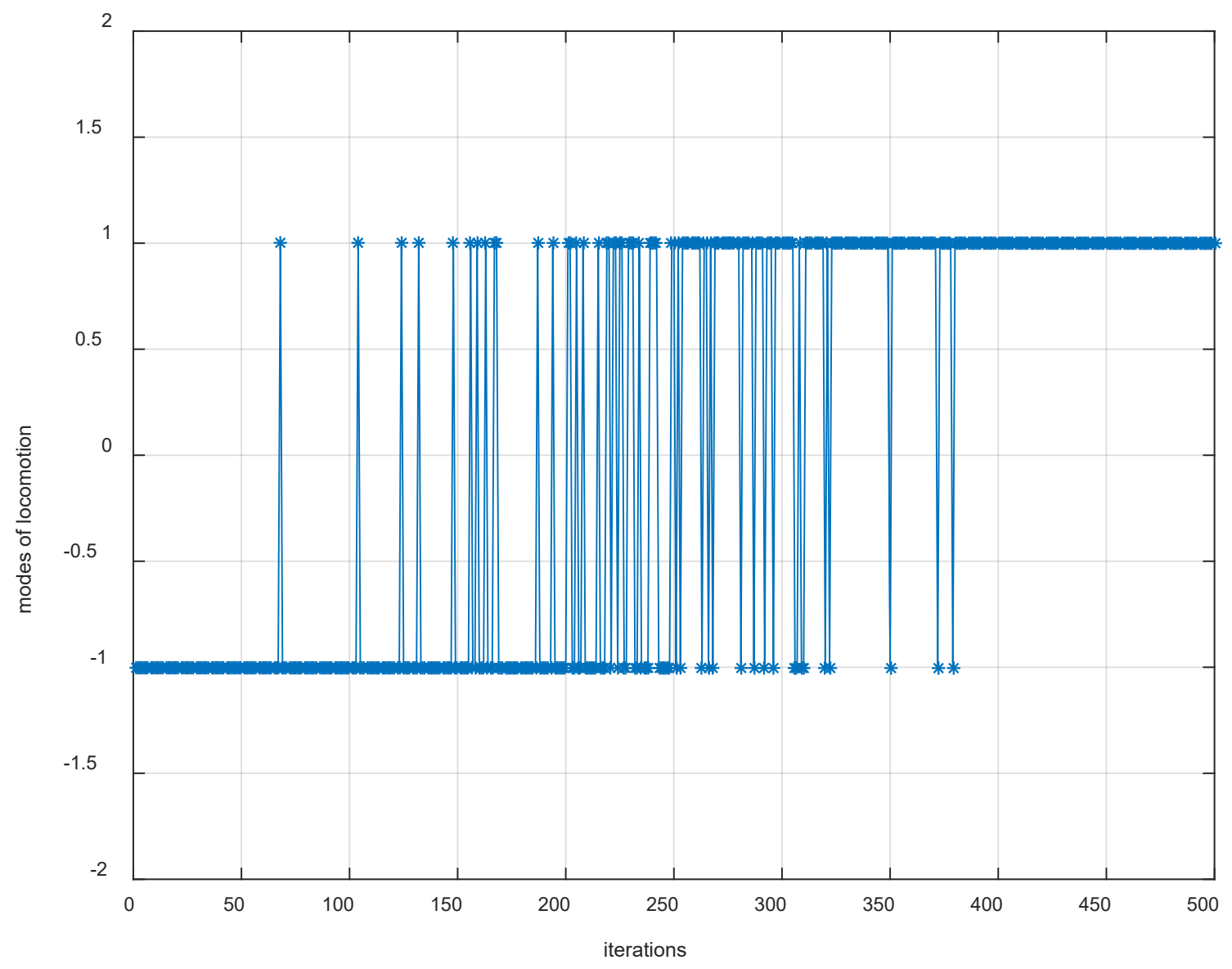

Fig. 11. The choice of modes of locomotion for a snake in the searching process.

Figure 12 illustrates the trajectory in three dimensions, focusing on the first three elements of the variable vector. Evidently, the snake roams the search space during the initial phase but gradually converges towards the theoretical global optimum.

Complementing Figure 12, Figure 13 presents a two-dimensional trace for the snake, specifically highlighting the first two elements of the variable vector. Notably, as the snake approaches the theoretical global optimum, there is a more concentrated distribution of touch spots. This observation serves to confirm, from one perspective, the efficacy of the SLLS's search mechanism.

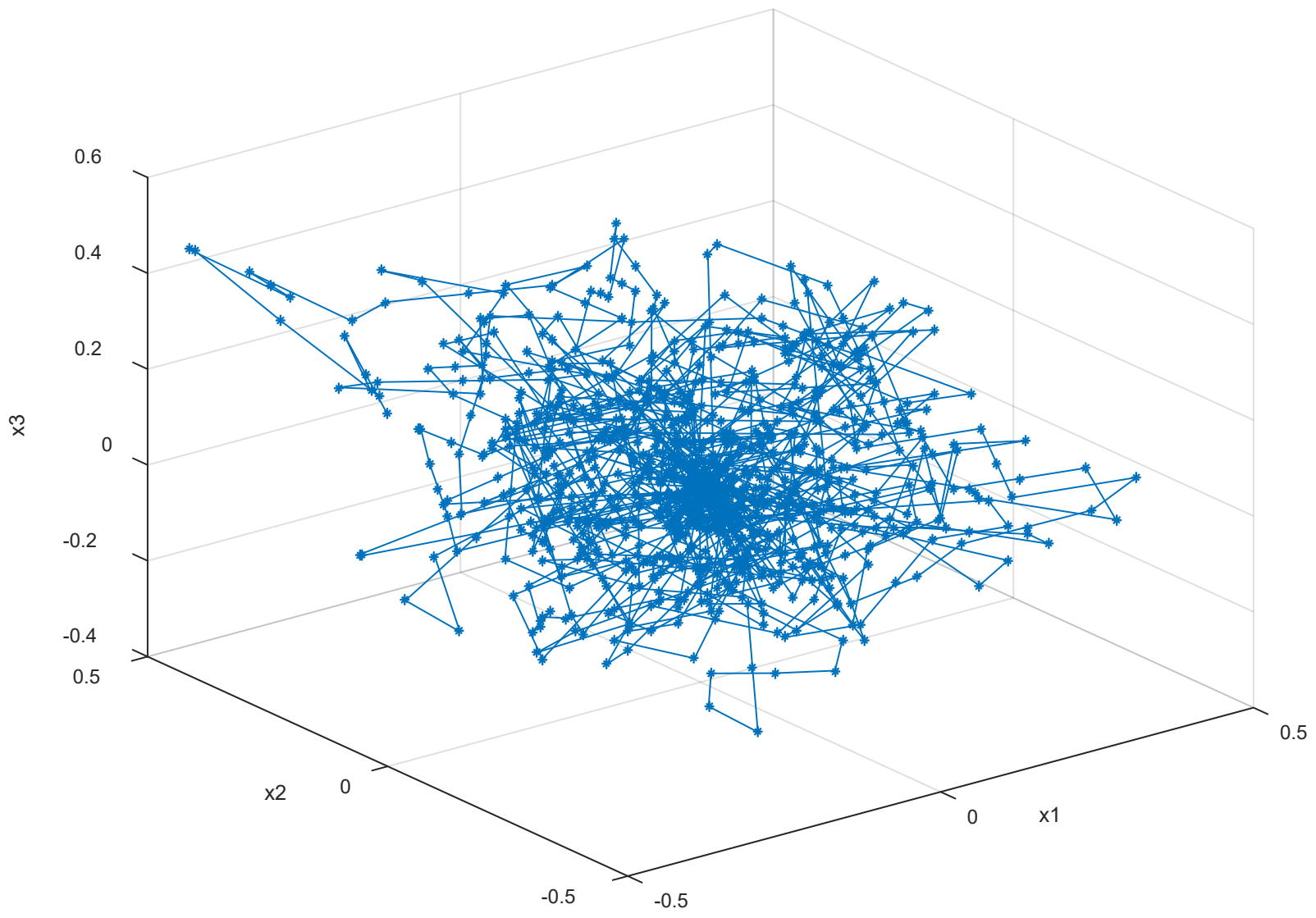

Fig. 12. The trace of the snake in 3D associated with its first three elements of the variable.

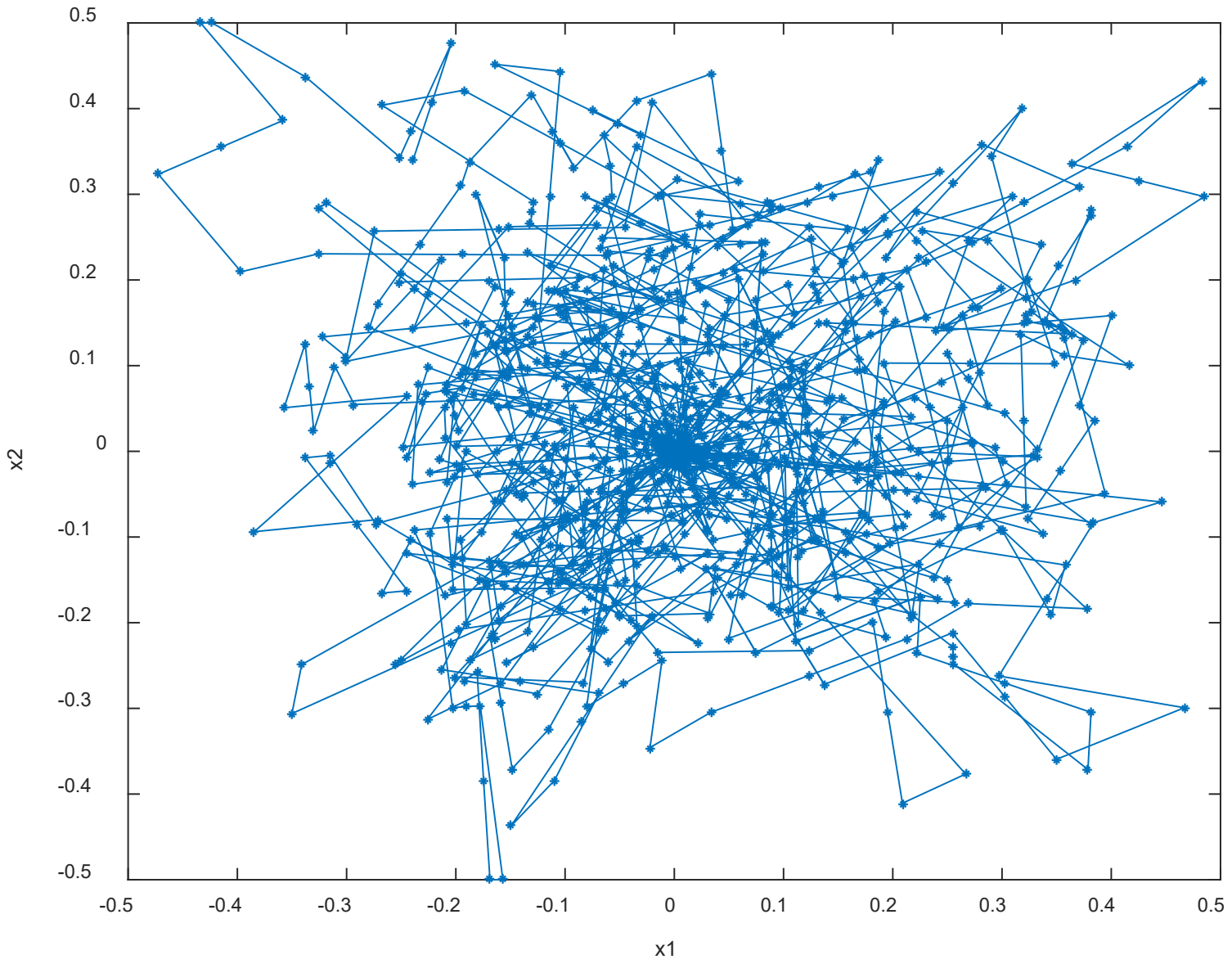

Fig. 13. The trace of the snake in 2D corresponding to its first two elements of the variable.

Let's commence the sensitivity analysis of the key parameters. Table 13 presents the results associated with different numbers of snakes in the swarm, ranging from 10 to 100. As the number of snakes increases, the average objective values generally decrease, but not in a strict linear fashion. Generally speaking, smaller numbers of snakes tend to yield larger best objective values compared to larger numbers of snakes. The worst objective values show a similar trend to the best objective values. It's evident that as the number of snakes increases, the Number of Total Movements (NTM) also increases accordingly. Given this observation, practitioners should exercise caution when attempting to improve final results by increasing the number of snakes, as this benefit comes at a higher computational cost.

**Table 13**

The performance of SLLS with changing number of snakes

| | Number of snakes | | | | | | | | | |
|---|---|---|---|---|---|---|---|---|---|---|
| | **10** | **20** | **30** | **40** | **50** | **60** | **70** | **80** | **90** | **100** |
| Mean | 0.001608 | 0.001284 | 7.861E-4 | 7.466E-4 | 7.363E-4 | 5.302E-4 | 5.955E-4 | 5.763E-4 | 5.965E-4 | 5.071E-4 |
| Std. | 6.678E-4 | 6.456E-4 | 5.549E-4 | 5.482E-4 | 4.395E-4 | 3.173E-4 | 4.759E-4 | 3.428E-4 | 3.648E-4 | 3.962E-4 |
| Best | 2.199E-4 | 2.387E-4 | 8.482E-5 | 1.612E-4 | 1.089E-4 | 3.342E-5 | 1.412E-5 | 8.557E-5 | 4.776E-5 | 1.260E-5 |
| Worst | 0.002937 | 0.002649 | 0.002451 | 0.002508 | 0.001940 | 0.001517 | 0.002457 | 0.001741 | 0.001696 | 0.001958 |
| NTM | 5000 | 10000 | 15000 | 20000 | 25000 | 30000 | 35000 | 40000 | 45000 | 50000 |

Table 14 presents the results associated with various radii for the truncated range of the natural growth function. As the radii increase from 5 to 50, the average objective values initially decrease, then rise later on. The best objective values exhibit a roughly similar trend. The worst objective values follow a pattern of decreasing as the radii increase from 5 to 25, after which they start to rise non-linearly with further increases in the radii. It's worth noting that when the radius takes on values of 30, 35, 40, or 45, the best objective value reaches 0.0, which is the optimal value for the problem.

The aforementioned observations suggest that selecting an appropriate value for the radius of the truncated range of the natural growth function leads to improved performance of the SLLS. This is because the radius adjusts the lengths of different learning phases. A larger radius results in longer starting and ending phases of the learning process and a shorter middle steep acceleration phase. The choice of a suitable radius enables the snakes to adopt movements at an appropriate rate as the learning process progresses.

**Table 14**

The performance of SLLS with changing radius of the truncated range given NTM=10,000.

| | The radius of the truncated range $\gamma$ | | | | | | | | | |
|---|---|---|---|---|---|---|---|---|---|---|
| | **5** | **10** | **15** | **20** | **25** | **30** | **35** | **40** | **45** | **50** |
| Mean | 0.001891 | 1.668E-4 | 8.881E-6 | 4.991E-9 | 5.109E-13 | 3.389E-7 | 1.042E-5 | 1.483E-5 | 1.070E-4 | 1.198E-4 |
| Std. | 8.998E-4 | 7.842E-5 | 4.268E-6 | 2.959E-9 | 4.471E-13 | 1.367E-6 | 4.346E-5 | 4.478E-5 | 2.361E-4 | 4.412E-4 |
| Best | 3.211E-4 | 6.492E-5 | 2.192E-6 | 4.268E-10 | 5.684E-14 | 0.0 | 0.0 | 0.0 | 0.0 | 5.329E-15 |
| Worst | 0.005120 | 3.595E-4 | 1.805E-5 | 1.305E-8 | 1.968E-12 | 7.059E-6 | 2.415E-4 | 2.401E-4 | 0.001220 | 0.002468 |

Table 15 illustrates the impact of the maximal number of iterations ($T$) on the performance of the SLLS. An increase in $T$ implies a higher computational cost. As $T$ increases, the average objective values generally decrease in a non-strict manner. The worst objective values follow a similar trend. However, unlike the average and worst objective values, the best objective values do not exhibit a consistent trend with increased $T$. The smallest of the best objective values corresponds to the condition of $T$=300. This observation suggests that achieving better performance with the SLLS does not entail simply increasing the maximal number of iterations, particularly when such an increase results in higher computational costs.

**Table 15**

The performance of SLLS with changing maximal number of iterations of SLLS.

| | The maximal number of iterations $T$ | | | | | | | | | |
|---|---|---|---|---|---|---|---|---|---|---|
| | **100** | **200** | **300** | **400** | **500** | **600** | **700** | **800** | **900** | **1000** |
| Mean | 0.003879 | 0.002308 | 0.001470 | 0.001274 | 0.001194 | 0.001041 | 0.001055 | 8.122E-4 | 9.900E-4 | 8.684E-4 |
| Std. | 0.001629 | 0.001204 | 8.184E-4 | 7.174E-4 | 6.830E-4 | 5.590E-4 | 6.308E-4 | 4.803E-4 | 6.622E-4 | 5.559E-4 |
| Best | 0.001027 | 4.278E-4 | 3.566E-6 | 7.970E-5 | 2.011E-4 | 2.064E-4 | 9.250E-5 | 1.039E-4 | 3.797E-5 | 1.734E-5 |
| Worst | 0.007901 | 0.005101 | 0.003978 | 0.002657 | 0.002916 | 0.002431 | 0.002527 | 0.001874 | 0.002706 | 0.002034 |
| NTM | 2000 | 4000 | 6000 | 8000 | 10000 | 12000 | 14000 | 16000 | 18000 | 20000 |

The list of visible spots stores historical global information, and its size directly influences the utility of this historical global information. Table 16 illustrates the impact of the size of the list of visible spots on the performance of the SLLS. The data reveals that when the size is not too small, such as 3 or 4, the average objective values are consistently relatively small, with no clear increasing or decreasing trend. The worst objective values associated with different sizes of the list of visible spots exhibit a similar pattern to the average objective values.

Interestingly, the SLLS achieves its best objective value when the size of the list equals 4. In all other cases, except when the size equals 4, the best objective values for larger sizes tend to be smaller than those for smaller sizes. This observation suggests that when implementing the SLLS, choosing a relatively large list size is a viable option. This choice can lead to relatively good performance without significantly increasing computational costs.

**Table 16**

The performance of SLLS with changing size of the list of visible points $m$.

| | **The size of the list of visible points $m$** | | | | | | | | | |
|---|---|---|---|---|---|---|---|---|---|---|
| | **3** | **4** | **5** | **6** | **7** | **8** | **9** | **10** | **11** | **12** |
| Mean | 0.001915 | 0.001385 | 9.721E-4 | 9.782E-4 | 0.001082 | 9.262E-4 | 8.392E-4 | 9.087E-4 | 9.419E-4 | 7.193E-4 |
| Std. | 8.324E-4 | 8.748E-4 | 4.047E-4 | 5.002E-4 | 4.944E-4 | 4.740E-4 | 5.867E-4 | 4.754E-4 | 4.583E-4 | 4.155E-4 |
| Best | 4.733E-4 | 3.522E-6 | 1.454E-4 | 1.779E-4 | 8.304E-5 | 1.270E-4 | 6.388E-5 | 1.760E-4 | 8.060E-5 | 8.285E-5 |
| Worst | 0.004086 | 0.002937 | 0.001895 | 0.002010 | 0.002278 | 0.001980 | 0.002445 | 0.002085 | 0.001788 | 0.001673 |

Table 17 illustrates the impacts of the minimal amplitude, $L_A^{min}$, and the number of half circles in a serpentine movement on the performance of the SLLS. Two significant observations are worth highlighting here.

Firstly, it's evident that the value of $L_A^{min}$ should be sufficiently small to ensure better SLLS performance. A comparison between the results with $L_A^{min}$ set to 1.0E-3 and the other two values of $L_A^{min}$ clearly shows that the condition $L_A^{min}$=1.0E-3 leads to notably worse results.

Secondly, for a given $L_A^{min}$, the results associated with different numbers of half circles in a serpentine movement do not exhibit a consistent increasing or decreasing trend as the number of half circles increases.

The reasons behind these observations are as follows. Regarding the first observation, a small $L_A^{min}$ makes the serpentine movement more effective during the final phase of the search process. In this phase, snakes should invest more energy in exploiting the more promising search areas. Concerning the second observation, the number of half circles in a serpentine movement determines how many touch points will be checked during the movement. It's important to note that an increase in touch points in a serpentine movement also translates to higher computational costs for a given maximal number of iterations. Similar to the analysis regarding the maximal number of iterations, the increase in the number of half circles has an indistinct impact on the final results. The variation in results without a clear trend simply reflects the random nature of the algorithm.

Based on the first observation, when setting the initial values for $L_A^{min}$, it is advisable to keep it sufficiently small. Regarding the second observation, it is recommended to experiment with different values for the number of half circles in serpentine locomotion if improved performance is desired.

**Table 17**

The performance of SLLS with the changed number of half circles in one serpentine locomotion given NTM=10,000.

| $L_A^{min}$ | | The number of half circles in one serpentine locomotion | | | | | | | | | |
|---|---|---|---|---|---|---|---|---|---|---|---|
| | | **1** | **2** | **3** | **4** | **5** | **6** | **7** | **8** | **9** | **10** |
| 1.0E-30 | Mean | 0.001121 | 0.001196 | 0.001146 | 0.001087 | 0.001273 | 0.001130 | 0.001202 | 0.001162 | 0.001564 | 0.001274 |
| | Std. | 5.369E-4 | 6.259E-4 | 6.102E-4 | 6.067E-4 | 7.990E-4 | 6.089E-4 | 7.724E-4 | 7.7275E-4 | 8.096E-4 | 6.788E-4 |
| | Best | 1.013E-4 | 3.312E-4 | 2.616E-5 | 2.040E-4 | 2.582E-4 | 9.622E-5 | 7.406E-5 | 1.516E-4 | 1.823E-4 | 8.987E-5 |
| | Worst | 0.002282 | 0.002603 | 0.002998 | 0.002920 | 0.003565 | 0.002636 | 0.003383 | 0.003265 | 0.002961 | 0.003042 |
| 1.0E-15 | Mean | 0.001228 | 0.001347 | 0.001103 | 0.001054 | 0.001157 | 0.001071 | 0.001378 | 0.001299 | 0.001023 | 0.001125 |
| | Std. | 5.365E-4 | 7.753E-4 | 5.908E-4 | 6.029E-4 | 5.907E-4 | 6.581E-4 | 7.693E-4 | 6.658E-4 | 5.248E-4 | 5.601E-4 |
| | Best | 2.489E-4 | 9.432E-5 | 1.332E-4 | 9.125E-5 | 1.505E-4 | 8.404E-5 | 2.119E-4 | 2.003E-4 | 1.814E-4 | 3.956E-5 |
| | Worst | 0.003106 | 0.003277 | 0.002505 | 0.002543 | 0.002570 | 0.002523 | 0.002945 | 0.003216 | 0.002400 | 0.002637 |
| 1.0E-3 | Mean | 0.002636 | 0.002733 | 0.002582 | 0.003320 | 0.002888 | 0.003323 | 0.002585 | 0.003669 | 0.003261 | 0.003044 |
| | Std. | 0.001637 | 0.001981 | 0.001503 | 0.001718 | 0.001556 | 0.001910 | 0.001484 | 0.002300 | 0.001872 | 0.002360 |
| | Best | 2.778E-4 | 1.024E-4 | 1.263E-4 | 0.001141 | 5.762E-4 | 4.122E-4 | 1.579E-5 | 5.894E-4 | 3.387E-4 | 3.268E-4 |
| | Worst | 0.005822 | 0.009618 | 0.006060 | 0.009567 | 0.005967 | 0.008193 | 0.006179 | 0.010897 | 0.008625 | 0.01028 |

Table 18 presents the impact of the number of touch points in a caterpillar locomotion on the performance of the SLLS. The data on average objective values associated with different numbers of touch points in a caterpillar movement reveal that both excessively small and excessively large numbers do not yield better objective values. Remarkably, when the number equals 5, the average objective value is minimized. The trend in worst objective values mirrors that of the average values.

In contrast to the patterns observed in average and worst objective values, the best objective values associated with different numbers of touch points in a caterpillar movement exhibit irregular changes. Despite the fact that excessively small or excessively large numbers of touch points still correspond to relatively high objective values, the variation in best objective values does not follow a clear trend. Based on this observation, it is advisable to determine an appropriate value for the number of touch points in a caterpillar locomotion through trial and error.

**Table 18**

The performance of SLLS with changing number of touch points in one caterpillar locomotion given NTM=10,000.

| | The number of touch points in one caterpillar locomotion | | | | | | | | | |
|---|---|---|---|---|---|---|---|---|---|---|
| | **3** | **4** | **5** | **6** | **7** | **8** | **9** | **10** | **11** | **12** |
| Mean | 0.001534 | 0.001090 | 8.509E-4 | 9.335E-4 | 0.001056 | 0.001162 | 0.001007 | 0.001157 | 0.001270 | 0.001748 |
| Std. | 6.668E-4 | 6.622E-4 | 4.728E-4 | 6.561E-4 | 8.809E-4 | 9.584E-4 | 6.992E-4 | 7.640E-4 | 9.404E-4 | 7.583E-4 |
| Best | 2.528E-4 | 1.025E-4 | 1.305E-4 | 1.600E-5 | 4.439E-5 | 5.389E-5 | 3.195E-5 | 1.042E-4 | 5.637E-8 | 2.980E-4 |
| Worst | 0.003058 | 0.003397 | 0.002006 | 0.002283 | 0.003763 | 0.003744 | 0.003112 | 0.003056 | 0.003632 | 0.003141 |

Finally, we will investigate the impact of the rate for demarcation, denoted as $r_{CL}$, on the performance of the SLLS. As illustrated earlier, when the value of $r_{CL}$ increases, all the touch points in a caterpillar movement tend to be closer to the target spot. Notably, when $r_{CL}$ is close to 1, the touch points in a caterpillar movement become very near to the target spot. Essentially, a caterpillar movement with a value of $r_{CL}$ near 1 resembles a neighborhood search with the target spot as the center. This can be advantageous if the target spot is close to the problem's optimum. However, it may diminish the exploratory potential of a caterpillar movement, especially when the search process has not reached its concluding phase. The data in Table 19 validate the aforementioned analysis.

Based on these observations, when determining an appropriate value for $r_{CL}$, opting for a value close to 0.5 can be a judicious choice.

**Table 19**

The performance of SLLS with changing rate for demarcation $r_{CL}$ given NTM=10,000.

| | The rate for demarcation $r_{CL}$ | | | | | | |
|---|---|---|---|---|---|---|---|
| | **0.3** | **0.4** | **0.5** | **0.6** | **0.7** | **0.8** | **0.9** |
| Mean | 0.001408 | 0.001273 | 0.001153 | 0.001147 | 0.001277 | 0.001463 | 0.001838 |
| Std. | 5.989E-4 | 6.535E-4 | 6.394E-4 | 5.331E-4 | 7.821E-4 | 7.527E-4 | 7.960E-4 |
| Best | 2.982E-4 | 2.401E-4 | 7.903E-5 | 1.578E-4 | 6.599E-5 | 1.178E-4 | 4.370E-4 |
| Worst | 0.002536 | 0.002792 | 0.002682 | 0.002531 | 0.002961 | 0.003102 | 0.003819 |

In the following, we will summarize the main observations from this subsection. The impacts of parameters can be broadly categorized into two sets.

In the first set, an optimal value for any parameter within this group typically falls in the middle range of feasible values. In other words, an ideal value for a parameter in this set should neither be excessively small nor too large. This first set comprises the radius of the truncated range of the natural growth function, the number of touch points in the caterpillar locomotion, the rate for demarcation, and the size of the list of visible points. Among these four parameters, the first three exhibit a consistent pattern, wherein the performance of the SLLS generally improves initially and then declines as the value of the respective parameter increases. However, the last parameter, specifically the number of half circles in a serpentine locomotion, does not strictly adhere to this pattern.

The remaining parameters, which include the number of half circles in a serpentine locomotion, the number of snakes in the swarm, and the maximal number of iterations, belong to the second set. Increasing the values of these parameters does not significantly enhance the performance of the SLLS, even though it comes at an added computational cost. When determining their values, practitioners should carefully strike a balance between the subtle improvement in final results and the associated computational expenses.

Lastly, it is crucial to ensure that the value of the minimal amplitude $L_A^{min}$ is sufficiently small, as a larger value of the minimal amplitude typically results in inferior outcomes compared to a smaller one.

## 5. **Conclusions**

In this paper, we have introduced a novel meta-heuristic algorithm known as the Snake Locomotion Learning Search algorithm (SLLS), inspired by the behavioral patterns observed in snakes. By simulating two distinct modes of snake locomotion, we have designed the exploration and exploitation mechanisms within this innovative algorithm. Serpentine locomotion primarily serves to enhance the exploration capabilities of the SLLS, while caterpillar locomotion, augmented by a trail-marking mechanism utilizing pheromones, is harnessed to maximize exploitation. To ensure a balanced exploration-exploitation dynamic throughout the search process, we employ a dynamically changing learning efficiency, generated by a natural growth rule, to determine the appropriate locomotion mode for each iteration of the SLLS.

The effectiveness and efficiency of the SLLS have been conclusively demonstrated through the resolution of 60 benchmark optimization problems and seven well-established engineering optimization challenges. In the majority of cases, the SLLS has surpassed alternative algorithms, as evidenced by superior final best and average objective values. In accordance with the "No Free Lunch Theorem" [2], these results underscore the theoretical and practical significance of the SLLS, promising its potential to address a wide array of optimization problems.

Moreover, we have conducted a meticulous analysis of the main parameters influencing the SLLS's final outcomes. Our observations have informed discussions on how to judiciously determine optimal parameter values and, in some instances, elucidated the reasons behind the observed impacts of parameter adjustments.

The adaptable design of the SLLS positions it for application across diverse optimization problems in the future. It holds promise in tackling complex tasks, such as real-size structural optimization. Furthermore, we envisage potential enhancements in the SLLS through the incorporation of additional snake locomotion modes, notably sidewinding locomotion and concertina locomotion, in future iterations. Additionally, a comprehensive study to delineate the operational mechanisms of the SLLS relative to other existing heuristic algorithms can provide valuable insights within a broader context.

## Appendix A (29 Standard Benchmark Optimization Problems):

**Table 1**

Unimodal benchmark functions.

| Functions | Dim | Range | $f_{\min}$ |
|---|---|---|---|
| $f_1(x) = \sum_{i=1}^{n} x_i^2$ | 30 | $[-100,100]^n$ | $f(0,0,\ldots,0) = 0$ |
| $f_2(x) = \sum_{i=1}^{n} \lvert x_i \rvert + \prod_{i=1}^{n} \lvert x_i \rvert$ | 30 | $[-10,10]^n$ | $f(0,0,\ldots,0) = 0$ |
| $f_3(x) = \sum_{i=1}^{n} (\sum_{j=1}^{i} x_j)^2$ | 30 | $[-100,100]^n$ | $f(0,0,\ldots,0) = 0$ |
| $f_4(x) = max\{\lvert x_i \rvert, 1 \le i \le n\}$ | 30 | $[-100,100]^n$ | $f(0,0,\ldots,0) = 0$ |
| $f_5(x) = \sum_{i=1}^{n-1} [100(x_{i+1} - x_i^2)^2 + (x_i - 1)^2]$ | 30 | $[-30,30]^n$ | $f(1,1,\ldots,1) = 0$ |
| $f_6(x) = \sum_{i=1}^{n} ([x_i + 0.5])^2$ | 30 | $[-100,100]^n$ | 0 |
| $f_7(x) = \sum_{i=1}^{n} i x_i^4 + random[0,1)$ | 30 | $[-1.28,1.28]^n$ | 0 |

**Table 2**

Multimodal benchmark functions.

| Functions | Dim | Range | $f_{\min}$ |
|---|---|---|---|
| $f_8(x) = \sum_{i=1}^{n} \left[-x_i \sin\left(\sqrt{\lvert x_i \rvert}\right)\right]$ | 30 | $[-500,500]^n$ | $-418.9829n$ |
| $f_9(x) = \sum_{i=1}^{n} [x_i^2 - 10\cos(2\pi x_i) + 10]$ | 30 | $[-5.12,5.12]^n$ | 0 |
| $f_{10}(x) = -20\exp\left(-0.2\sqrt{\frac{1}{n}\sum_{i=1}^{n} x_i^2}\right) - \exp\left(\frac{1}{n}\sum_{i=1}^{n}\cos(2\pi x_i)\right) + 20 + e$ | 30 | $[-32,32]^n$ | 0 |
| $f_{11}(x) = \sum_{i=1}^{n} \frac{x_i^2}{4000} - \prod_{i=1}^{n} \cos\left(\frac{x_i}{\sqrt{i}}\right) + 1$ | 30 | $[-512,512]^n$ | 0 |
| $f_{12}(x) = \frac{\pi}{n}\{10\sin^2(\pi y_1) + \sum_{i=1}^{n-1}(y_i - 1)^2[1 + 10\sin^2(\pi y_{i+1})] + (y_n - 1)^2\} + \sum_{i=1}^{n} u(x_i, 10, 100, 4)$ | 30 | $[-50,50]^n$ | 0 |
| $y_i = 1 + \frac{x_i + 1}{4}$, $u(x_i, a, k, m) = \begin{cases} k(x_i - a)^m & x_i > a \\ 0 & -a \le x_i \le a \end{cases}$ | | | |
| $f_{13}(x) = 0.1\left\{\sin^2(3\pi x_1) + \sum_{i=1}^{n}(x_i - 1)^2[1 + \sin^2(3\pi x_i + 1)] + (x_n - 1)^2[1 + \sin^2(2\pi x_n)]\right\} + \sum_{i=1}^{n} u(x_i, 5, 100, 4)$ | 30 | $[-50,50]^n$ | 0 |

**Table 3**

Fixed dimension multimodal benchmark functions.

| Functions | Dim | Range | $f_{\min}$ |
|---|---|---|---|
| $f_{14}(x) = \left(\frac{1}{500} + \sum_{j=1}^{25} \frac{1}{j+\sum_{i=1}^{2}(x_i - a_{ij})^6}\right)^{-1}$ | 2 | $[-65.536,65.536]^n$ | 1 |
| $f_{15}(x) = \sum_{i=1}^{11} \left[a_i - \frac{x_1(b_i^2+b_i x_2)}{b_i^2+b_i x_3+x_4}\right]^2$ | 4 | $[-5,5]^n$ | 0.0003075 |
| $f_{16}(x) = 4x_1^2 - 2.1x_1^4 + x_1^6/3 + x_1x_2 - 4x_2^2 + 4x_2^4$ | 2 | $[-5,5]^n$ | -1.0316285 |
| $f_{17}(x) = \left(x_2 - \frac{5.1}{4\pi}x_1^2 + \frac{5}{\pi}x_1 - 6\right)^2 + 10\left(1 - \frac{1}{8\pi}\right)\cos x_1 + 10$ | 2 | $[-5,5]^n$ | 0.398 |
| $f_{18}(x) = [1 + (x_1 + x_2 + 1)^2(19 - 14x_1 + 3x_1^2 - 14x_1 + 6x_1x_2 + 3x_2^2)] \times$ $[30 + (2x_1 - 3x_2)^2(18 - 32x_1 + 12x_1^2 + 48x_2 - 36x_1x_2 + 27x_2^2)]$ | 2 | $[-2,2]^n$ | 3 |
| $f_{19}(x) = -\sum_{i=1}^{4} c_i \exp\left(-\sum_{j=1}^{3} a_{ij}(x_j - p_{ij})^2\right)$ | 3 | $[0,1]^n$ | -3.86 |
| $f_{20}(x) = -\sum_{i=1}^{4} c_i \exp\left(-\sum_{j=1}^{6} a_{ij}(x_j - p_{ij})^2\right)$ | 6 | $[0,1]^n$ | -3.32 |
| $f_{21}(x) = -\sum_{i=1}^{5}[(x - a_i)(x - a_i)^T + c_i]^{-1}$ | 4 | $[0,10]^n$ | -10.1532 |
| $f_{22}(x) = -\sum_{i=1}^{7}[(x - a_i)(x - a_i)^T + c_i]^{-1}$ | 4 | $[0,10]^n$ | -10.4028 |
| $f_{23}(x) = -\sum_{i=1}^{10}[(x - a_i)(x - a_i)^T + c_i]^{-1}$ | 4 | $[0,10]^n$ | -10.5363 |

**Table 4**

Composite benchmark functions.

| Functions | Dim | Range | $f_{\min}$ |
|---|---|---|---|
| $F_{24}$ (CF1)<br>$f_1, f_2, f_3, \ldots, f_{10}$ =Sphere Function<br>$[\sigma_1, \sigma_2, \sigma_3, \ldots, \sigma_{10}] = [1,1,1,\ldots,1]$<br>$[\lambda_1, \lambda_2, \lambda_3, \ldots, \lambda_{10}] = [5/100, 5/100, 5/100, \ldots, 5/100]$ | 10 | [-5,5] | 0 |
| $F_{25}$ (CF2)<br>$f_1, f_2, f_3, \ldots, f_{10}$ =Griewank's Function<br>$[\sigma_1, \sigma_2, \sigma_3, \ldots, \sigma_{10}] = [1,1,1,\ldots,1]$<br>$[\lambda_1, \lambda_2, \lambda_3, \ldots, \lambda_{10}] = [5/100, 5/100, 5/100, \ldots, 5/100]$ | 10 | [-5,5] | 0 |
| $F_{26}$ (CF3)<br>$f_1, f_2, f_3, \ldots, f_{10}$ =Griewank's Function<br>$[\sigma_1, \sigma_2, \sigma_3, \ldots, \sigma_{10}] = [1,1,1,\ldots,1]$<br>$[\lambda_1, \lambda_2, \lambda_3, \ldots, \lambda_{10}] = [1,1,1,\ldots,1]$ | 10 | [-5,5] | 0 |
| $F_{27}$ (CF4)<br>$f_1, f_2$= Ackley's Function, $f_3, f_4$= Rastrigin's Function,<br>$f_5, f_6$= Weierstrass Function, $f_7, f_8$= Griewank's Function,<br>$f_9, f_{10}$= Sphere Function<br>$[\sigma_1, \sigma_2, \sigma_3, \ldots, \sigma_{10}] = [1,1,1,\ldots,1]$<br>$[\lambda_1, \lambda_2, \lambda_3, \ldots, \lambda_{10}]$ =[5/32, 5/32, 1, 1, 5/0.5, 5/0.5, 5/100, 5/100, 5/100, 5/100] | 10 | [-5,5] | 0 |
| $F_{28}$ (CF5)<br>$f_1, f_2$= Rastrigin s Function, $f_3, f_4$= Weierstrass Function,<br>$f_5, f_6$= Griewank s Function, $f_7, f_8$= Ackley s Function,<br>$f_9, f_{10}$ = Sphere Function<br>$[\sigma_1, \sigma_2, \sigma_3, \ldots, \sigma_{10}] = [1,1,1,\ldots,1]$<br>$[\lambda_1, \lambda_2, \lambda_3, \ldots, \lambda_{10}]$ =[1/5, 1/5, 5/0.5, 5/0.5, 5/100, 5/100, 5/32, 5/32, 5/100, 5/100] | 10 | [-5,5] | 0 |
| $f_{29}$ (CF6)<br>$f_1, f_2$= Rastrigin's Function, $f_3, f_4$= Weierstrass' Function,<br>$f_5, f_6$= Griewank's Function, $f_7, f_8$= Ackley's Function,<br>$f_9, f_{10}$= Sphere Function<br>$[\sigma_1, \sigma_2, \sigma_3, \ldots, \sigma_{10}]$ =[0.1, 0.2, 0.3, 0.4, 0.5, 0.6, 0.7, 0.8, 0.9, 1]<br>$[\lambda_1, \lambda_2, \lambda_3, \ldots, \lambda_{10}]$ =[0.1 ∗ 1/5, 0.2 ∗ 1/5, 0.3 ∗ 5/0.5, 0.4 ∗ 5/0.5, 0.5 ∗ 5/100, 0.6 ∗ 5/100, 0.7 ∗ 5/32, 0.8 ∗ 5/32, 0.9 ∗ 5/100, 1 ∗ 5/100] | 10 | [-5,5] | 0 |

# Appendix B: (7 Classic Engineering Optimization Problems)

**B01. Multiple disc clutch brake design**

Consider variable $x = (r_i, r_o, t, F, Z)$

Minimize $f(x) = \pi(r_o^2 - r_i^2)t(Z+1)\rho$

Subject to:

$g_1(x) = r_o - r_i - \Delta r \geq 0,\ g_2(x) = l_{max} - (Z+1)(t+\delta) \geq 0,\ g_3(x) = p_{max} - p_{rz} \geq 0,$

$g_4(x) = p_{max}v_{srmax} - p_{rz}v_{sr} \geq 0,\ g_5(x) = v_{srmax} - v_{sr} \geq 0,\ g_6(x) = T_{max} - T \geq 0,$

$g_7(x) = M_h - sM_s \geq 0$, $g_8(x) = T \geq 0$,

where $M_h = \frac{2}{3}\mu F Z \frac{r_o^3 - r_i^3}{r_o^2 - r_i^2}$, $p_{rz} = \frac{F}{\pi(r_o^2 - r_i^2)}$, $M_h = \frac{2\pi n(r_o^3 - r_i^3)}{90(r_o^2 - r_i^2)}$, $T = \frac{I_z \pi n}{30(M_h + M_f)}$,

$\Delta r = 20$mm, $t_{max}$=3mm, $t_{min}$=1.5mm, $l_{max} = 30$mm, $Z_{max}$=10, $v_{srmax}$=10m/s, $\mu = 0.5$,

$s = 1.5$, $M_s$=40Nm, $M_f$=3Nm, n=250rpm, $p_{max}$=1MPa, $I_z$=55kg mm$^2$, $T_{max}$=15s, $F_{max}$=1000N,

$r_{min}$=55mm, $r_{o\,max}$=110mm, $\rho = 7800 \text{kg}/m^3$.

Variable range $r_i \in \{60, 61, 62, \ldots, 80\}$, $r_o \in \{90, 91, 92, \ldots, 110\}$, $t \in \{1.0, 1.5, 2.0, 2.5, 3\}$, $F \in \{600, 610, 620, \ldots, 1000\}$, $Z \in \{2,3,4,5,6,7,8,9\}$.

**B02. Robot gripper**

Minimize $f(x) = \max_z F_k(x, z) - \min_z F_k(x, z)$

Subject to:

$g_1(x) = Y_{min} - y(x, Z_{max}) \geq 0$, $g_2(x) = y(x, Z_{max}) \geq 0$, $g_3(x) = y(x, 0) - Y_{max} \geq 0$,

$g_4(x) = Y_G - y(x, 0) \geq 0$, $g_5(x) = (a + b)^2 - l^2 - e^2 \geq 0$,

$g_6(x) = (l - Z_{max})^2 + (a - e)^2 - b^2 \geq 0$, $g_7(x) = l - Z_{max} \geq 0$,

$g_8(x) = d(l, Z_{max}, e) + b - a \geq 0$,

$g_9(x) = d(l, 0, e) + b - a \geq 0$,

$g_{10}(x) = b + a - d(l, 0, e) \geq 0$,

where $d(l, z, e) = \sqrt{(l - z)^2 + e^2}$, $\emptyset = \tan^{-1}\left(\frac{e}{l - z}\right)$, $\alpha = \cos^{-1}\left(\frac{a^2 + d^2 - b^2}{2ad}\right) + \emptyset$,

$\beta = \cos^{-1}\left(\frac{b^2 + d^2 - a^2}{2bd}\right) - \emptyset$, $F_k = \frac{Pb\sin(\alpha + \beta)}{2c\cos\alpha}$, $y(x, z) = 2(e + f + c\,sin(\beta + \delta))$,

$Y_{min} = 50$, $Y_{max}$=100, $Y_G$=150, $Z_{max}$=100, $P = 100$,

$10 \leq a, b, f \leq 150$, $100 \leq c \leq 200$, $0 \leq e \leq 50$, $100 \leq l \leq 300$, $1 \leq \delta \leq 3.14$.

The constraints $g_8(x)$, $g_9(x)$ and $g_{10}(x)$ are newly added to guarantee the realizability of the results.

**B03. Rolling element bearing**

Maximize $C_d = f_c Z^{2/3} D_b^{1.8}$ if $D_b \leq 25.4$mm

$C_d = 3.64 f_c Z^{2/3} D_b^{1.4}$ if $D_b > 25.4$mm

Subject to:

$g_1(x) = \frac{\phi_o}{2\sin^{-1}(D_b/D_m)} - Z + 1 \geq 0$, $g_2(x) = 2D_b - K_{Dmin}(D - d) \geq 0$,

$g_3(x) = K_{Dmax}(D - d) - 2D_b \geq 0$, $g_4(x) = \zeta B_w - D_b \leq 0$,

$g_5(x) = D_m - 0.5(D + d) \geq 0$, $g_6(x) = (0.5 + e)(D + d) - D_m \geq 0$,

$g_7(x) = 0.5(D - D_m - D_b) - \varepsilon D_b \geq 0$, $g_8(x) = f_i \geq 0.515$, $g_9(x) = f_o \geq 0.515$,

where

$$f_c = 37.91\left\{1 + \left[1.04\left(\frac{1-\gamma}{1+\gamma}\right)^{1.72}\left(\frac{f_i(2f_o - 1)}{f_o(2f_i - 1)}\right)^{0.41}\right]^{10/3}\right\}^{-0.3}\left(\frac{\gamma^{0.3}(1-\gamma)^{1.39}}{(1+\gamma)^{1/3}}\right)\left(\frac{2f_i}{2f_i - 1}\right)^{0.41},$$

$$\phi_o = 2\pi - 2\cos^{-1}\frac{\left[\frac{(D-d)}{2} - 3\left(\frac{T}{4}\right)\right]^2 + \left(\frac{D}{2} - \frac{T}{4} - D_b\right)^2 - \left(\frac{d}{2} + \frac{T}{4}\right)^2}{2\left[\frac{(D-d)}{2} - 3\left(\frac{T}{4}\right)\right]\left(\frac{D}{2} - \frac{T}{4} - D_b\right)},$$

$\gamma = \frac{D_b\cos\alpha}{D_m}$, $f_i = \frac{r_i}{D_b}$, $f_o = \frac{r_o}{D_b}$, $T = D - d - 2D_b$, $D = 160$, $d$=90, $B_w$=30, $\alpha = 0$,

$0.5(D + d) \leq D_m \leq 0.6(D + d)$, $0.15(D - d) \leq D_b \leq 0.45(D - d)$, $Z \in \{4,5,6,\ldots,50\}$,

$0.515 \leq f_i, f_o \leq 0.6$, $0.4 \leq K_{Dmin} \leq 0.5$, $0.6 \leq K_{Dmax} \leq 0.7$, $0.3 \leq \varepsilon \leq 0.4$,

$0.02 \leq e \leq 0.1$, $0.6 \leq \zeta \leq 0.85$.

**B04. Hydrodynamic thrust bearing**

Minimize $f(x) = \frac{QP_0}{0.7} + E_f$

Subject to: $g_1(x) = W - W_s \geq 0$, $g_2(x) = P_{max} - P_0 \geq 0$, $g_3(x) = \Delta T_{max} - \Delta T \geq 0$, $g_4(x) = h - h_{min} \geq 0$, $g_5(x) = R - R_0 \geq 0$, $g_6(x)=0.001-\frac{\gamma}{gP_0}\left(\frac{Q}{2\pi Rh}\right)^2 \geq 0$, $g_7(x) = 5000 - \frac{W}{\pi(R^2-R_0^2)} \geq 0$,

Where $\gamma = 0.0307$, C=0.5, n=-3.55, $C_1$=10.04, $W_s$=101000, $P_{max}$=1000, $\Delta T_{max}$=50, $h_{min}$=0.001,
$g = 386.4$, N=750, $1 \leq R, R_0, Q \leq 16$, $1e-6 \leq \mu \leq 16e-6$.

**B05. Belleville spring**

Minimize $f(x) = 0.07075\pi(D_e^2 - D_i^2)t$

Subject to:

$g_1(x)=S - \frac{4E\delta_{max}}{(1-\mu^2)\alpha D_e^2}\left[\beta\left(h - \frac{\delta_{max}}{2}\right) + \gamma t\right] \geq 0$,
$g_2(x)=\frac{4E\delta_{max}}{(1-\mu^2)\alpha D_e^2}\left[\left(h - \frac{\delta_{max}}{2}\right)(h - \delta_{max})t + t^3\right] - P_{max} \geq 0$,
$g_3(x)=\delta_l - \delta_{max} \geq 0$, $g_4(x)=H - h - t \geq 0$, $g_5(x) = D_{max} - D_e \geq 0$,
$g_6(x) = D_e - D_i \geq 0$, $g_7(x) = 0.3 - \frac{1}{D_e - D_i} \geq 0$,

where $\alpha = \frac{6}{\pi \ln K}\left(\frac{K-1}{K}\right)^2$, $\beta = \frac{6}{\pi \ln K}\left(\frac{K-1}{\ln K} - 1\right)$, $\gamma = \frac{6}{\pi \ln K}\left(\frac{K-1}{2}\right)$,

$P_{max} = 5400$lb, $\delta_{max} = 0.2$in., $S$=200kPsi, $E$=30e6psi, $\mu = 0.3$, $H$=0.2in., $D_{max}$=12.01in.,
$K = \frac{D_e}{D_i}$, $\delta_l = f(a)h$, $a = h/t$. $0.01 \leq t \leq 6.0$, $0.05 \leq h \leq 0.5$, $5.0 \leq D_i, D_e \leq 15.0$.

Values of f(a) vary as shown in Table 7.

Table 7

Variation of $f(a)$ with $a$.

| $a$ | ≤ 1.4 | 1.5 | 1.6 | 1.7 | 1.8 | 1.9 | 2.0 | 2.1 | 2.2 | 2.3 | 2.4 | 2.5 | 2.6 | 2.7 | ≥ 2.8 |
|---|---|---|---|---|---|---|---|---|---|---|---|---|---|---|---|
| $f(a)$ | 1 | 0.85 | 0.77 | 0.71 | 0.66 | 0.63 | 0.6 | 0.58 | 0.56 | 0.55 | 0.53 | 0.52 | 0.51 | 0.51 | 0.5 |

**B06. Step-cone pulley**

Minimize $f(x) = \rho w\left\{d_1^2\left[1 + \left(\frac{N_1}{N}\right)^2\right] + d_2^2\left[1 + \left(\frac{N_2}{N}\right)^2\right] + d_3^2\left[1 + \left(\frac{N_3}{N}\right)^2\right] + d_4^2\left[1 + \left(\frac{N_4}{N}\right)^2\right]\right\}$

Subject to:

$h_1(x) = C_1 - C_2 = 0$, $h_2(x) = C_1 - C_3 = 0$, $h_3(x) = C_1 - C_4 = 0$,
$g_{1,2,3,4}(x) = R_i \geq 2$, $g_{5,6,7,8}(x) = P_i \geq (0.75 * 745.6998)$,
where $C_i$ indicates the length of the belt to obtain speed $N_i$ and is given by

$C_i = \frac{\pi d_i}{2}\left(1 + \frac{N_i}{N}\right) + \frac{\left(\frac{N_i}{N}-1\right)^2 d_i^2}{4a} + 2a$, i = (1,2,3,4).

$R_i$ is the tension ratio and is given by
$R_i = \exp(\mu\left[\pi - 2\sin^{-1}\left(\left(\frac{N_i}{N} - 1\right)\frac{d_i}{2a}\right)\right])$, i = (1,2,3,4).

$P_i$=stw$\left[1 - exp\left[-\mu\left\{\pi - 2\sin^{-1}\left(\left(\frac{N_i}{N} - 1\right)\frac{d_i}{2a}\right)\right\}\right]\right]\frac{\pi d_i N_i}{60}$, i = (1,2,3,4).

$\rho = 7200\text{kg}/m^3$, $a = 3$m, $\mu = 0.35$, s = 1.75MPa, t = 8mm.
$1\text{mm} \leq d_1, d_2, d_3, d_4, w \leq 100\text{mm}$.

**B07. Speed reducer design**

Minimize $f(x) = 0.7854x_1x_2^2(3.3333x_3^2 + 14.9334x_3 - 43.0934) - 1.508x_1(x_6^2 + x_7^2)$
$+7.4777(x_6^3 + x_7^3) + 0.7854(x_4x_6^2 + x_5x_7^2)$

Subject to:

$g_1(x) = \frac{27}{x_1x_2^2x_3} - 1 \leq 0$, $g_2(x) = \frac{397.5}{x_1x_2^2x_3^2} - 1 \leq 0$, $g_3(x) = \frac{1.93x_4^3}{x_2x_3x_6^4} - 1 \leq 0$, $g_4(x) = \frac{1.93x_5^3}{x_2x_3x_7^4} - 1 \leq 0$,

$g_5(x) = \frac{\sqrt{[745x_4/(x_2x_3)]^2+16.9e6}}{110x_6^3} - 1 \leq 0,\ \ g_6(x) = \frac{\sqrt{[745x_5/(x_2x_3)]^2+157.5e6}}{85x_7^3} - 1 \leq 0,$

$g_7(x) = \frac{x_2x_3}{40} - 1 \leq 0,\ g_8(x) = \frac{5x_2}{x_1} - 1 \leq 0,\ g_9(x) = \frac{x_1}{12x_2} - 1 \leq 0,$

$g_{10}(x) = \frac{1.5x_6+1.9}{x_4} - 1 \leq 0,\ g_{11}(x) = \frac{1.1x_7+1.9}{x_5} - 1 \leq 0,$

where $2.6 \leq x_1 \leq 3.6;\ 0.7 \leq x_2 \leq 0.8;\ 17 \leq x_3 \leq 28;\ 7.3 \leq x_4 \leq 8.3;\ 7.8 \leq x_5 \leq 8.3;$

$2.9 \leq x_6 \leq 3.9;\ 5.0 \leq x_7 \leq 5.5.$